\documentclass[11pt, reqno]{amsart}
\RequirePackage[OT1]{fontenc}
\usepackage[colorlinks=true, urlcolor=blue, linkcolor=blue, citecolor=blue, pdftex]{hyperref}
\usepackage{csquotes}
\topmargin - 1 cm
\usepackage[utf8]{inputenc}
\usepackage[matrix,arrow,curve]{xy}
\usepackage{amsmath,amsfonts,verbatim,afterpage,euscript,mathrsfs,amssymb}
\usepackage{array}
\usepackage{amsthm}
\usepackage{dsfont}
\usepackage{hyperref}
\usepackage{enumerate}
\usepackage[english]{babel}
\usepackage{tikz}
\usetikzlibrary{shapes}
\usepackage{titlesec}
\usepackage{algorithmicx}
\usepackage{algpseudocode}
\usepackage{algorithm}

\titleformat{\section}[block]{\normalfont\bfseries\filcenter}{\itshape\thesection}{1em}{}
\titleformat{\subsection}[block]{\normalfont\bfseries}{\itshape\thesubsection}{0.9em}{}
\titleformat{\subsubsection}[block]{\normalfont\bfseries}{\itshape\thesubsubsection}{0.8em}{}
\titleformat{\caption}[block]{\normalfont}{\itshape}{0.8em}{}

\newcommand{\mc}[1]{\ensuremath{\mathcal{#1}}}

\newcommand{\tX}{\tilde{X}}
\newcommand{\tq}{\tilde{q}}

\newcommand{\E}{\mathbb{E}}

\renewcommand{\P}{\mathbb{P}}
\newcommand{\R}{\mathbb{R}}

\renewcommand{\L}{\mathbb{L}}

\newcommand{\mL}{\mathcal{L}}

\newcommand{\mE}{\mathcal{E}}
\newcommand{\mP}{\mathcal{P}}

\newcommand{\mI}{\mathcal{I}}
\newcommand{\mH}{\mathcal{H}}
\newcommand{\bbu}{\mathbf{u}}
\newcommand{\p}{\partial}

\renewcommand{\mE}{\mathcal{E}}
\renewcommand{\p}{\partial}

\newcommand{\ind}{\mathbf{1}_{\mathcal{S}}}
\newcommand{\indc}{\mathbf{1}_{\mathcal{S}^c}}
\renewcommand{\d}{{\rm d}}

\numberwithin{equation}{section}

\newtheorem{theoreme}{Theorem}[section]
\newtheorem{proposition}[theoreme]{Proposition}

\newtheorem{lemme}[theoreme]{Lemma}
\newtheorem{definition}[theoreme]{Definition}
\newtheorem*{remarque}{Remark}

\newtheorem{claim}[theoreme]{Claim}
\begin{document}
\title[Strong existence and uniqueness for degenerate SDE with H\"{o}lder drift]{Strong existence and uniqueness for stochastic differential equation with H\"{o}lder drift and degenerate noise}

\author{P.E. Chaudru de Raynal}
\address{Université Nice Sophia Antipolis, CNRS, LJAD, UMR 7351, 06100 Nice, France.}
\email[P.E. Chaudru de Raynal]{pe.deraynal@univ-smb.fr}
\date {September, 2015}

\begin{abstract}
In this paper, we prove pathwise uniqueness for stochastic degenerate systems with a H\"{o}lder drift, for a H\"{o}lder exponent larger than the critical value $2/3$. This work extends to the degenerate setting the earlier results obtained by Zvonkin \cite{zvonkin_transformation_1974}, Veretennikov \cite{veretennikov_strong_1980}, Krylov and R\"{o}ckner \cite{krylov_strong_2005} from non-degenerate to degenerate cases. The existence of a threshold for the H\"{o}lder exponent in the degenerate case may be understood as the price to pay to balance the degeneracy of the noise. Our proof relies on regularization properties of the associated PDE, which is degenerate in the current framework and is based on a parametrix method.
\end{abstract}

\maketitle

\section{Introduction}
Let $T$ be a positive number and  $d$ be a positive integer, we consider the following $\mathbb{R}^d \times \mathbb{R}^d$ system for any $t$ in $[0,T]$:

\begin{equation}\label{systemEDO}
\left\lbrace \begin{array}{llll}
\d X^1_t = F_1(t,X^1_t,X^2_t)\d t + \sigma(t,X^1_t,X_t^2) \d W_t,\qquad & X^1_0=x_1,\\
\d X^2_t = F_2(t,X^1_t,X^2_t)\d t,\qquad & X^2_0=x_2,\\
\end{array}
\right.
\end{equation}
where $x_1$, $x_2$ belong to $\R^d$, $(W_{t}, t\geq 0)$ is a standard $d$-dimensional Brownian motion defined on some filtered probability space $(\Omega, \mathcal{F}, \mathbb{P}, (\mathcal{F}_t)_{ t \geq 0})$ and $F_1,F_2,\sigma: [0,T] \times \mathbb{R}^d \times \mathbb{R}^d  \to \mathbb{R}^d \times \mathbb{R}^d \times \mathcal{M}_{d}(\mathbb{R})$ (the set of real $d \times d$ matrices) are measurable functions. The diffusion matrix $a:=\sigma \sigma^*$ is supposed to be uniformly elliptic. The notation ``$*$'' stands for the transpose.\\ 

In this paper, we investigate the well posedness of \eqref{systemEDO} outside the Cauchy-Lipschitz framework. Notably, we are  interested in the strong posedness, \emph{i.e.} strong existence and uniqueness of a solution. Strong existence means that there exists a process $(X^1_t,X^2_t,\ 0\leq t \leq T)$ adapted to the filtration generated by the Brownian motion $(W_{t},\ 0\leq t \leq T)$ which satisfies \eqref{systemEDO}. Strong uniqueness means that if two processes satisfy this equation with the same initial conditions, their trajectories are almost surely indistinguishable. Here, we show that under a suitable H\"{o}lder assumption on the drift coefficients and Lipschitz condition on the diffusion matrix, the strong well-posedness holds for \eqref{systemEDO}.\\

It may be a real challenge to prove the existence of a unique solution for a differential system without Lipschitz conditions on the coefficients. For example, in \cite{zhang_well-posedness_2010} Zhang showed that SDE with Sobolev coefficients admits a unique generalized solution (as Lebesgue almost everywhere stochastic flow, see Definition 2.1 in \cite{zhang_well-posedness_2010}) under integrability conditions on the drift, the divergence of the drift and the diffusion coefficient. The case of drift with divergence of polynomial growth being also handled.\

Concerning strong solutions, the first result in this direction is due to Zvonkin. In \cite{zvonkin_transformation_1974}, the author showed that the strong well-posedness holds for  the one-dimensional system 
\begin{equation}\label{simpleEq}
Y_t = \int_0^t b(s,Y_s)\d s + W_t,\ Y_0=y \quad t\in [0,T],
\end{equation} 
for a measurable function $b$ in $\L^{\infty}$. Then, Veretennikov \cite{veretennikov_strong_1980} generalized this result to the multi-dimensional case. Krylov and R\"{o}ckner showed in \cite{krylov_strong_2005} the strong well-posedness for $b$ in $\L^{p}_{{\rm loc}},\ p>d$ and Zhang in \cite{zhang_stochastic_2011} handled the case of non-constant, Sobolev and non-degenerate diffusion coefficient. Finally, when $b$ is a measurable and bounded function, Davie showed in \cite{davie_uniqueness_2007} that for almost every Brownian path, there exists a unique solution for the system \eqref{simpleEq}. We emphasize that this result implies the strong uniqueness, but the converse is not true. Indeed, in such a case, there exists an \emph{a priori} set $\Omega' \subset \Omega$ with $\P(\Omega')=1$ such that for all $\omega$ in $\Omega'$ the solution of \eqref{simpleEq} is unique.\\
 
All these results rely on the regularization of differential systems by adding a non-degenerate noise, and we refer to \cite{flandoli_random_2011} for a partial review on this subject. The proofs of such results rely on the deep connection between SDEs and PDEs (see \cite{bass_diffusions_1998} or \cite{friedman_stochastic_2006} for a partial review in the elliptic and parabolic cases). The generator associated to the Markov process $Y$ is a linear partial differential operator of second order (usually denoted by $\mL$) with the transition density of $Y$ as fundamental solution. As explained by Fedrizzi and Flandoli in \cite{fedrizzi_pathwise_2011}: ``if we have a good theory for the PDE:
\begin{equation}\label{pdez}
\frac{\p}{\p t}u + \mL u = \Phi,\ \text{on }[0,T] \quad u_T=\mathbf{0},
\end{equation}
where the source term $\Phi$ has the same regularity as the drift, then, we have the main tools to prove strong uniqueness''.\\

 In \eqref{systemEDO}, the noise added is completely degenerate w.r.t the component $X^2$. This sort of system has also been studied by Veretennikov in \cite{veretennikov_stochastic_1983} but without considering any regularization in the degenerate direction. Indeed, the author showed that strong well-posedness holds when the drift is measurable and bounded and the diffusion matrix is Lipschitz w.r.t the non-degenerate component $X^1$ and when both the drift and the diffusion matrix are twice continuously differentiable functions with bounded derivatives w.r.t the degenerate component.\\

In this paper, we show that the noise regularizes, even in the degenerate direction, by means of the random drift. Unfortunately, there is a price to pay to balance the degeneracy of the noise. First, the drift must be at least $2/3$-H\"{o}lder continuous w.r.t the degenerate component. We do not know how sharp is this critical value, but it is consistent with our approach. Secondly, the drift $F_2$ of the second component must be Lipschitz continuous w.r.t the first component and its derivative in this direction has to be uniformly non degenerate: this allows the noise to propagate through the system and then the drift to regularize.\\

Our proof also relies on regularization properties of the associated PDE, and the aforementioned ``good theory'' is here a ``strong theory'': a Lipschitz bound on the solution of \eqref{pdez} and on its derivative w.r.t the first component. We emphasize that, in our case, the generator $\mL$ is given by: for all $\psi$ in $C^{1,2,1}([0,T] \times \mathbb{R}^{d} \times \mathbb{R}^{d},\R^d)$\footnote{\emph{i.e. continuously differentiable w.r.t. the first variable, twice continuously differentiable w.r.t. the second variable and once continuously differentiable w.r.t. the third variable.}}
\begin{eqnarray}
\mL\psi(t,x_1,x_2) &=& \frac{1}{2}{\rm Tr}(a(t,x_1,x_2)D^2_{x_1} \psi(t,x_1,x_2))  + \left[F_1(t,x_1,x_2)\right] \cdot \left[D_{x_1}\psi (t,x_1,x_2)\right]\nonumber\\
&& \qquad + \left[F_2(t,x_1,x_2)\right] \cdot \left[D_{x_2} \psi (t,x_1,x_2)\right] .\label{gengen}
\end{eqnarray} 
where ${\rm Tr}(a)$ stands for the trace of the matrix $a$, ``$\cdot$'' denotes the standard Euclidean inner product on $\R^d$ and where for any $z$ in $\R^d$, the notation $D_{z}$ means the derivative w.r.t the variable $z$. Here, the crucial point is that the operator is not uniformly parabolic. 

When the coefficients are smooth and when the Lie algebra generated by the vector fields spans the whole space, it is known that such an operator admits a smooth fundamental solution (see \cite{hormander_hypoelliptic_1967}): it is said to be hypoelliptic and the coefficients are said to satisfy a H\"{o}rmander condition. The assumption on the uniform non-degeneracy of the derivative of the drift $F_2$ w.r.t $x_1$ together with the uniform ellipticity of the matrix $\sigma$ can be understood as a sort of weak H\"{o}rmander condition.\\

This particular form of degeneracy is a non-linear generalization of Kolmogorov's degeneracy, in reference to the first work \cite{kolmogorov_zufallige_1934} of Kolmogorov in this direction. Degenerate operators of this form have been studied by many authors see e.g. the work of Di Francesco and Polidoro \cite{di_francesco_schauder_2006}, and Delarue and Menozzi \cite{delarue_density_2010}. We also emphasize that, in \cite{menozzi_parametrix_2011}, Menozzi deduced  from the regularization property exhibited in (\cite{delarue_density_2010}) the well weak posedness of a generalization of \eqref{systemEDO}. Nevertheless, to the best of our knowledge, there does not exist a strong theory,  in the sense defined above, for the PDE \eqref{pdez} when $\mL$ is defined by \eqref{gengen}. We investigate it by using the so called parametrix approach (see \cite{friedman_partial_1964} for partial review in the elliptic setting).\\

To conclude this introduction, we just come back to the regularity assumed on the drift of the non-degenerate component. In comparison with the works of Veretennikov \cite{veretennikov_strong_1980,veretennikov_stochastic_1983}, Krylov and R\"{o}ckner \cite{krylov_strong_2005}, and Flandoli and Fedrizzi \cite{fedrizzi_pathwise_2011}, asking for $F_1$ to be only in $\L^p$, $p>d$ might appear as the right framework. Since the parametrix is a perturbation method and since we are interested in $\mathbb{L}^{\infty}$ estimates, we suppose the drift $F_1$ to be H\"{o}lder continuous w.r.t $x_1$.\\

\subsection{Organization of this paper}
Subsection \ref{NOT} states useful notations. In Subsection \ref{GSAMR} we give the detailed assumptions and the main result of this paper: strong existence and uniqueness for \eqref{systemEDO}. In Subsection \ref{SOP}, we expose the strategy to prove this result. It is based on the regularization properties of the associated PDE which are, in fact, the main contribution of this work. These properties are given in Subsection \ref{SecPDEres}. Finally, our main result is proved in Subsection \ref{POT}.\\ 
The remainder of this paper is dedicated to the proof of the regularization properties of the associated PDE.

We present in Section \ref{LBH} the linear and Brownian heuristic. It explains how the proof of the regularization properties in a simple case works and allows to understand our assumptions and how the proof in the general case can be achieved. Then, we give in Section \ref{SOTLSOP} the mathematical tools, and the proof of the regularization properties of the PDE is given in Section \ref{PBFGH}. This is the technical part of this paper.

\subsection{Notations}\label{NOT} In order to simplify the notations, we adopt the following convention: $x, y, z,\xi,$ \emph{etc.} denote the $2d-$dimensional real variables $(x_1,x_2), (y_1,y_2), (z_1,z_2), (\xi_1,\xi_2),$ \emph{etc.}.  Consequently, each component of the $d$-dimensional variables $x_k,\ k=1,2$ is denoted by $x_{kl},\ l=1,\cdots,d$. We denote by $g(t,X_t)$ any function $g(t,X_t^1,X_t^2)$ from $[0,T] \times\R^{d} \times \R^d$ to $\R^{N},\ N\in \mathbb{N}$. Here, $X_t = (X_t^1,X_t^2)$ and then $F(t,X_t)$ is the $\R^{2d}$ valued function $(F_1(t,X_t^1,X_t^2), F_2(t,X_t^1,X_t^2))^*$. We rewrite the system (\ref{systemEDO}) in a shortened form:
\begin{equation}\label{simplisys}
\d X_t = F(t,X_t)\d t  + B\sigma(t,X_t)\d W_t,
\end{equation}
where $B$ is the $2d \times d$ matrix: $B=(\rm{Id}, 0_{\R^{d}\times \R^d})^*$. ``${\rm Id}$'' stands for the identity matrix of $\mathcal{M}_d(\R)$, the set of real $d \times d$ matrices. When necessary, we write $(X^{t,x}_s)_{t\leq s \leq T}$ the process in \eqref{systemEDO} which starts from $x$ at time $t$, \emph{i.e.} such that $X_t^{t,x}=x$.

We recall that the canonical Euclidean inner product on $\R^d$ is denoted by ``$\cdot$''. We denote by ${\rm GL}_d(\R)$ the set of $d\times d$ invertible matrices with real coefficients and by $\phi$ a measurable function from $[0,T] \times \R^d \times \R^{d}$ to $\R^2$. Each one-dimensional component of this function is denoted by $\phi_i$, $i=1,2$ and plays the role of one coordinate of $F_i$. Hence, $\phi_i$ satisfies the same regularity as $F_i$ given latter. We recall that $a$ denotes the square of the diffusion matrix $\sigma$, $a :=\sigma \sigma^*$. Subsequently, we denote by $c$, $C$, $c'$, $C'$, $c''$ \emph{etc.} a positive constant, depending only on known parameters in \textbf{(H)}, given just below, that may change from line to line and from an equation to another.

We denote by $C^{1,2,1}$ the space of functions that are continuously differentiable w.r.t. the first variable, twice continuously differentiable w.r.t. the second variable and and once continuously differentiable w.r.t. the third variable.

The notation $D$ stands for the total space derivative. For any function from $[0,T]\times \R^d \times \R^d$ we denote by $D_1$ (resp. $D_2$) the derivative with respect to the first (resp. second) $d$-dimensional space component. In the same spirit, the notation $D_{z}$ means the derivative w.r.t the variable $z$. Hence, for all integer $n$, $D^n_{z}$ is the $n^{\rm{th}}$ derivative w.r.t $z$ and for all integer $m$ the $n\times m$ cross differentiations w.r.t $z$, $y$ are denoted by $D^n_{z}D^m_{y}$. Furthermore, the partial derivative $\p/\p_t$ is denoted by $\p_t$.\\

\subsection{Main Result}\label{GSAMR}
\textbf{Assumptions} \textbf{(H).} We say that assumptions \textbf{(H)} hold if the following assumptions are satisfied:
\begin{description}
\item[(H1) regularity of the coefficients] There exist $0<\beta_i^j<1$, $1\leq i,j\leq2$ and three positive constants $C_1, C_2, C_\sigma$ such that for all $t$ in $[0,T]$ and all $(x_1,x_2)$ and $(y_1,y_2)$ in $\mathbb{R}^d \times \mathbb{R}^d$
\begin{eqnarray*} 
&&|F_1(t,x_1,x_2) - F_1(t,y_1,y_2)| \leq C_{1} (|x_1-y_1|^{\beta^1_1} + |x_2-y_2|^{\beta^2_1})\\
&&|F_2(t,x_1,x_2) - F_2(t,y_1,y_2)| \leq C_2(|x_1-y_1|+|x_2-y_2|^{\beta^2_2})\\
&&|\sigma(t,x_1,x_2) - \sigma(t,y_1,y_2)| \leq C_\sigma (|x_1-y_1| + |x_2-y_2|).
\end{eqnarray*}
Moreover, the coefficients are supposed to be continuous w.r.t the time and the exponents $\beta_i^2,\ i=1,2$ are supposed to be strictly greater than $2/3$. Thereafter, we set $\beta_2^1=1$ for notational convenience.
\item[(H2) uniform ellipticity of $\sigma\sigma^*$] The function $\sigma\sigma^*$ satisfies the uniform ellipticity hypothesis:
\begin{equation*}
\exists \Lambda >1,\ \forall \zeta \in \mathbb{R}^{2d},\quad \Lambda^{-1}|\zeta|^2\leq \left[\sigma \sigma^*(t,x_1,x_2)\zeta\right] \cdot \zeta  \leq \Lambda |\zeta|^2,
\end{equation*}
for all $(t,x_1,x_2) \in  [0,T] \times \mathbb{R}^d \times \mathbb{R}^d$.
\item[(H3-a) differentiability and regularity of  $\ x_1 \mapsto F_2(.,x_1,.)$]  For all $(t,x_2) \in  [0,T]\times \R^d$, the function $F_2(t,.,x_2):\ x_1 \mapsto F_2(t,x_1,x_2)$ is continuously differentiable and there exist $0<\alpha^1<1$ and a positive constant $ \bar{C}_2$ such that, for all $(t,x_2)$ in $[0,T] \times \mathbb{R}^d$ and $x_1,y_1$ in $\mathbb{R}^d$
\begin{eqnarray*} 
&&|D_{1}F_2(t,x_1,x_2) - D_{1}F_2(t,y_1,x_2)| \leq \bar{C}_2 |x_1-y_1|^{\alpha^1}.
\end{eqnarray*} 
\item[(H3-b) non degeneracy of $(D_{1}F_2)(D_{1}F_2)^*$] There exists a closed convex subset $\mathcal{E} \subset {\rm GL}_d(\R)$  (the set of $d\times d$ invertible matrices with real coefficients) such that for all $t$ in $[0,T]$ and $(x_1,x_2)$ in $\R^{2d}$ the matrix $D_{1}F_2(t,x_1,x_2)$ belongs to $\mathcal{E}$. We emphasize that this implies that 
\begin{equation*}
\exists \bar{\Lambda} >1,\ \forall \zeta \in \mathbb{R}^{2d},\quad \bar{\Lambda}^{-1}|\zeta|^2\leq \left[(D_{1}F_2)(D_{1}F_2)^*(t,x_1,x_2)\zeta\right] \cdot \zeta  \leq \bar{\Lambda} |\zeta|^2,
\end{equation*}
for all $(t,x_1,x_2) \in  [0,T] \times \mathbb{R}^d \times \mathbb{R}^d$.
\end{description}

\begin{remarque} The reason for the existence of the critical value $2/3$ for the H\"{o}lder regularity of the drift in \textbf{(H1)} and the particular ``convexity'' assumption \textbf{(H3-b)} are discussed in Section \ref{LBH}. In the following, the sentence ``known parameters in \textbf{(H)}'' refers to the parameters belonging to these assumptions. 
\end{remarque}

The following Theorem is the main result of this paper and regards the strong well-posedness of the system \eqref{systemEDO}.

\begin{theoreme}\label{c1:MR}
Under \textbf{(H)}, strong existence and uniqueness hold for \eqref{systemEDO} for any positif $T$.
\end{theoreme}

\subsection{Strategy of proof}\label{SOP}
Let us expose the basic arguments to prove Theorem \ref{c1:MR}. Existence of a weak solution follows from a compactness argument. Then, if the strong uniqueness holds, the strong existence follows. The main issue consists in proving the strong uniqueness. 

This works as follows: suppose that there exists a unique $C^{1,2,1}([0,T]\times \R^d\times \R^d,\R^d)$ solution $\bbu=(\bbu_1,\bbu_2)^*$ of the  linear system of PDEs:
\begin{equation}\label{PDE1}
\left\lbrace
\begin{array}{lll}
\partial_t\bbu^{}_i(t,x) + \mL^{}\bbu^{}_i(t,x) = F^{}_i(t,x), \quad \text{for }(t,x) \in [0,T) \times \mathbb{R}^{2d},\\
\bbu^{}_i(T,x)=0_{\R^{d}},\quad i=1,2,
\end{array}
\right.
\end{equation}
then, thanks to It\^{o}'s formula, for all $t$ in $[0,T]$ we have that
$$\int_0^t F(s,X^{x}_s) \d s = \bbu(t,X^x_t) - \bbu(0,x) - \int_0^t D_{}\bbu(s,X^{x}_s)B\sigma(s,X^x_s) \d W_s.$$
Thus, if  the functions $\bbu$ and $D \bbu B = (D_{1}\bbu,0_{\R^d\times \R^d})^*$ are $C_T$ Lipschitz continuous, the drift $\int_0^t F(s,X^x_s) \d s$ is Lipschitz continuous w.r.t. the argument $X^x$, with Lipschitz constant $C_T$.

Now suppose that $C_T$ tends to 0 when $T$ goes to 0. Uniqueness then follows by mean of classical and circular arguments for $T$ small enough. Since the strategy can be iterated, one can deduce strong uniqueness on any positive interval.\\

The main issue here is then to obtain a strong theory for the PDE \eqref{PDE1}. The problem is that, to the best of our knowledge, it is not known that this PDE admits a Lipschitz $C^{1,2,1}([0,T]\times \R^d\times \R^d,\R^d)$ solution under our weak H\"{o}rmander assumptions. Nevertheless, in our analysis we do not need to obtain the existence of a regular solution, but only the existence of the Lipschitz bounds for $u$ and $D_{1}u$ depending only on known parameters in \textbf{(H)}. 

Therefore, we investigate these bounds in a smooth setting. Thanks to  assumptions \textbf{(H)}, there exists a sequence of smooth (say infinitely differentiable with bounded derivatives of all order greater than 1) mollified coefficients $(a^n,F_1^n,F_2^n)_{n \geq 1}$ satisfying \textbf{(H)} uniformly (in $n$), such that this sequence converges in supremum norm to $(a,F_1,F_2)$. More details on  the regularization procedure are given in Subsection \ref{SecPDEres} below. In this regularized setting, we know that the PDE \eqref{PDE1} admits a unique smooth solution.

Then, thanks to a first order parametrix expansion of the solution $\bbu^n$ of the regularized PDE \eqref{PDE1} (see \emph{e.g.} \cite{friedman_partial_1964}), we show that for $T$ small enough there exists a positive constant $C_T$, which is independent of the regularization procedure, such that the supremum norm of $D_{1}\bbu^n$, $D_{2}\bbu^n$, $D_{1}^2\bbu^n$ and $D_{2}D_{1}\bbu^n$ are bounded by $C_T$. We then recover the Lipschitz regularity of the drift in small time and we obtain strong uniqueness by letting the regularization procedure tends to the infinity.\\

\subsection{PDE's results}\label{SecPDEres}
This section summarizes the PDE's results used for proving Theorem \ref{c1:MR}.

\textbf{The mollifying procedure.} Let us first detail how the smooth approximation of the coefficients $a,F_1,F_2$ works. For all positive integer $n$, we set:
$$F_2^n(t,x) = \int_{} F_2(t-s,x-y) \varphi_1^n(y)\varphi_2^n(s)\d y \d s,$$
where $\varphi_1^n(.) = c_1n^{2d}\varphi(n|.|)$ and  $\varphi_2^n(.) = c_2n\varphi(n|.|)$ for $c_1$, $c_2$ two constants of normalization and for a smooth function $\varphi$ with support in the unit ball. For example $\varphi: z\in\R \mapsto \exp\left(-\frac{1}{1-z^2}\right)\mathbf{1}_{]-1;1[}(z)$. 

By defining $(F_1^n)_{n\geq 1}$ and $(a^n)_{n \geq 1}$ with the same procedure, it is then clear that for every $n$ the mollified coefficients $a^n,F_1^n,F_2^n$ are infinitely differentiable with bounded derivatives of all order greater than 1 and such that
\begin{equation}\label{defmoll}
(a^n,F_1^n,F_2^n) \underset{n\to+\infty}{\longrightarrow} (a,F_1,F_2),
\end{equation} 
uniformly on $[0,T] \times \R^d \times \R^d$. Moreover, it is well-seen that they satisfy the same assumptions as $(a,F_1,F_2)$ uniformly in $n$. 

Let us just check the non-degeneracy assumption $\textbf{(H3-b)}$ on $D_{1}F_2^n$. For all positive $\delta$, one can find a positive integer $N(\delta)$ and a sequence of rectangles $(R_k)_{1\leq k \leq N(\delta)}$ having sides of length less than $\delta$ and a family of points $\{(s_k,y_k)\in R_k,\ 1\leq k\leq N(\delta)\}$ such that, for all $(t,x)$ in $[0,T]\times \R^{2d}$:
$$D_{1}F_2^n(t,x) =\lim_{\delta \to 0} \sum_{k=1}^{N(\delta)} D_{1}F_{2}(t-s_k,x-y_k) \int_{R_k}\varphi_1^n(y)\varphi_2^n(s)\d y\d s.$$
Since $D_{1}F_2$ belongs to the closed convex subset $\mE$, it is clear that $D_{1}F_2^n$ belongs to $\mE$.\\

\textbf{The regularized PDE.}
As we said, we do not solve the limit PDE problem \eqref{PDE1}. The investigations are done with the mollified coefficients $(a^n,F_1^n,F_2^n)_{ n \geq 1}$ defined above. Let us denote by $\mL^n$ the regularized version of $\mL$ (that is the version of $\mL$ with mollified coefficients). We have from Section 2.1 Chapter II of \cite{freidlin_functional_1985} (note that the time dependence here is not a problem to do so):
\begin{lemme}\label{eeu}
Let $n$ be a positive integer. The PDE,
\begin{equation}\label{regPED}
\partial_t \bbu^n_i(t,x) + \mL^n\bbu^n_i(t,x) = F^n_i(t,x), \quad \text{for }(t,x) \in [0,T] \times \mathbb{R}^{2d} \quad \bbu^n_i(T,x)=0_{\R^d},\quad i=1,2,
\end{equation}
where  $\mL^n$ is the regularized version of the operator $\mL$ defined by \eqref{gengen}, admits a unique solution $\bbu^n=(\bbu_1^n,\bbu_2^n)^*$, which is infinitely differentiable. 
\end{lemme}

 Besides, the terminal condition $\bbu^n(T,.)=0_{\R^{2d}}$ is very important: it guarantees that the solution and its derivatives vanish at time $T$. Hence, it allows to control the Lipschitz constant of $\bbu^n$ by a constant small as $T$ is small. Indeed, we show in Section \ref{PBFGH} that the solutions $\bbu^n,\ n \geq 1$ satisfy:
\begin{proposition}\label{EE}
There exist a positive $\mc{T}$, a positive $\delta_{\ref{EE}}$ and a positive constant $C$, depending only on known parameters in \textbf{(H)} and not on $n$, such that, for all positive $T$ less than $\mc{T}$:
\begin{equation*}
||D_{1}\bbu^n||_{\infty} + ||D_{2}\bbu^n||_{\infty} + ||D^2_{1}\bbu^n||_{\infty} + ||D_{1}D_{2}\bbu^n||_{\infty} \leq CT^{\delta_{\ref{EE}}}.
\end{equation*}
\end{proposition}

In order to prove these results,  we emphasize that each coordinate of the vectorial solution $\bbu^n_i$ of the decoupled linear PDE (\ref{regPED}) can be described by the PDE
\begin{equation}\label{PDE2}
\partial_t u^n_i(t,x) + \mL^n u_i(t,x) = \phi^n_i(t,x), \quad \text{for }(t,x) \in [0,T] \times \mathbb{R}^{2d},\quad u^n_i(T,x)=0,\ i=1,2,
\end{equation}
where $\phi_i^n:\R^{2d}\to \R$ denotes the mollified (by the procedure described above) coefficient $\phi_i$: a function that satisfies the same regularity assumptions as $F_i$ given in \textbf{(H1)} (this function plays the role of one of the coordinates of $F_i$). Therefore, we only have to prove Lemma \ref{eeu} and Proposition \ref{EE} for \eqref{PDE2} instead of (\ref{regPED}).\\

Since the estimates on the solutions $\bbu^n,\ n \geq 1$ are obtained uniformly in $n$ (that is independently of the procedure of regularization), when we investigate the properties of the solution of the PDE \eqref{PDE2} in the following sections, we forget the superscript ``$n$'' which arises from the mollifying procedure, and we further assume that the following assumptions hold:

\textbf{Assumptions} \textbf{(HR).} We say that assumptions \textbf{(HR)} hold if assumptions \textbf{(H)} hold true and $F_1, F_2, \phi_1, \phi_2$ and $ a$ are infinitely differentiable functions with bounded derivatives of all order greater than 1. 

\subsection{Proof of Theorem \ref{c1:MR}}\label{POT}
We know from Theorem 6.1.7 of \cite{stroock_multidimensional_1979} that the system (\ref{systemEDO}) admits a weak solution (we emphasize that this result remains valid under the linear growth conditions assumed on the coefficients). Hence, we only have to prove the strong uniqueness. Thereafter, we denote by ``$\mathbf{1}$'' the $2d\times 2d$ matrix:
\begin{equation}\label{c1defdu1}
\begin{pmatrix}
{\rm Id} & 0_{\R^d \times \R^d}\\
0_{\R^d \times \R^d} & 0_{\R^d \times \R^d}
\end{pmatrix}.
\end{equation}

Let $(X_t, t\geq0)$ and $(Y_t, t\geq0)$ be two solutions of (\ref{systemEDO}) with the same initial condition $x$ in $\R^{2d}$. Let $\bbu^{n}$ be the solution of the linear system of PDEs (\ref{regPED}). Thanks to Lemma \ref{eeu}, we can apply It\^{o}'s formula on $\bbu^{n}(t,X_t)-X_t$ and we obtain

\begin{eqnarray*}
\bbu^n(t,X_t)-X_t &=&  \int_{0}^t \left[\p_t \bbu^n  + \mL^{}\bbu^n\right](s,X_s) \d s - \int_{0}^t F(s,X_s) \d s +\bbu^n(0,x)-x  \nonumber\\
&&\quad + \int_{0}^t \left[D_{}\bbu^n- \mathbf{1}\right]B\sigma (s,X_s) \d W_s.
\end{eqnarray*}
In order to use the fact that $\p_t \bbu^n  + \mL^{n}\bbu^n = F^n$, we rewrite

\begin{eqnarray*}
\bbu^n(t,X_t)-X_t &=&  \int_{0}^t \left[\p_t \bbu^n  + \mL^{n}\bbu^n\right](s,X_s) \d s + \int_{0}^t \left(\mL-\mL^{n}\right)\bbu^n(s,X_s) \d s\\
&&     - \int_{0}^t F(s,X_s) \d s +\bbu^n(0,x)-x + \int_{0}^t \left[D_{}\bbu^n- \mathbf{1}\right]B\sigma (s,X_s) \d W_s,
\end{eqnarray*}
and then,

\begin{eqnarray*}
\bbu^n(t,X_t)-X_t &=&  \int_{0}^t \left(\mL-\mL^{n}\right)\bbu^n(s,X_s) \d s   + \int_{0}^t (F^n(s,X_s)- F(s,X_s)) \d s +\bbu^n(0,x)-x  \nonumber\\
&&\quad + \int_{0}^t \left[D_{}\bbu^n- \mathbf{1}\right]B\sigma (s,X_s) \d W_s.
\end{eqnarray*}

By the same arguments, we have:
\begin{eqnarray*}
\bbu^n(t,Y_t)-Y_t &=&  \int_{0}^t \left(\mL-\mL^{n}\right)\bbu^n(s,Y_s) \d s  + \int_{0}^t (F^n(s,Y_s) -  F(s,Y_s)) \d s + \bbu^n(0,x)-x  \nonumber\\
&&\quad + \int_{0}^t \left[D_{}\bbu^n- \mathbf{1}\right]B\sigma (s,Y_s) \d W_s.
\end{eqnarray*}

By taking the expectation of the supremum over $t$ of the square norm of the difference of the two equalities above, we get from Doob's inequality that

\begin{eqnarray*}
\E \left[\sup_{t \in [0,T]}|X_t-Y_t|^2 \right]  &\leq & C\Bigg\{ \E\left[ \sup_{t \in [0,T]} | \bbu^n(t,X_t)-\bbu^n(t,Y_t)|^2\right]\\
&&  + \E \left[\int_0^T \left|\left[D_{}\bbu^nB\right](s,X_s)-\left[D_{}\bbu^nB\right](s,Y_s)\right|^2 \left\|\sigma\right\|_{\infty}^2 \d s \right]\\
&&  +  \E \left[\int_0^T (\left\| D_{}\bbu^nB \right\|_{\infty} + \mathbf{1}) \left|\left[\sigma(s,Y_s) - \sigma(s,X_s)\right]\right|^2 \d s \right] +  \mathcal{R}(n,T)\Bigg\},
\end{eqnarray*}
where
\begin{eqnarray*}
\mathcal{R}(n,T) &=&  T\Bigg( \E\left[\sup_{t \in [0,T]}  |F^{(n)}(t,Y_t)-F(t,Y_t)|^2\right] + \E\left[\sup_{t \in [0,T]}  |(\mL^{n}-\mL)\bbu^n(t,Y_t)|^2\right] \\
&& \quad +   \E\left[\sup_{t \in [0,T]}  |F^{(n)}(t,X_t)-F(t,X_t)|^2 \right] + \E\left[\sup_{t \in [0,T]}  |(\mL^{n}-\mL)\bbu^n(t,X_t)|^2 \right]\Bigg).
\end{eqnarray*}

First, note that from \eqref{defmoll}, for both $Y_t$ and $X_t$, we have:
$$\E\left[\sup_{t \in [0,T]}  |F^{n}(t,X_t)-F(t,X_t)|^2\right] + \E\left[\sup_{t \in [0,T]}  |(\mL^{n}-\mL)\bbu^n(t,X_t)|^2\right]   \to 0,\ \text{as }n \to \infty,$$
so that $\mathcal{R}(n,T)$ tends to $0$ when $n$ tends to $+\infty$.
Secondly, we know from Proposition \ref{EE}, that there exists a positive $\mc{T}$ such that for all $T$ less than $\mc{T}$ and for all $t$ in $ [0,T]$, the functions $\bbu^n$ and $D_{1}\bbu^n$ are Lipschitz continuous in space, with a Lipschitz constant independent of $n$. Since $D_{}\bbu^nB = (D_{1}\bbu^n,0_{\R^d\times \R^d})$, by letting $n$ tends to  $+\infty$ and using the two arguments above, we deduce that for all $T$ less than $\mc{T}$:
 
\begin{equation*}
\E \left[\sup_{t \in [0,T]}|X_t-Y_t|^2\right] \leq C(T) \left\lbrace  \E \left[\sup_{t \in [0,T]}|X_t-Y_t|^2\right] +\E \left[\int_0^T |X_s-Y_s|^2 \d s \right]\right\rbrace,
\end{equation*}
where $C(T)$ tends to $0$ when $T$ tends to $0$. Hence, one can now find another positive $\mc{T}'$, depending only on known parameters in \textbf{(H)}, such that strong uniqueness holds for all $T$ less than $\mc{T}'$. By iterating this computation, the same result holds on any finite intervals and so on $[0,\infty)$.\qed

\section{The linear and Brownian heuristic}\label{LBH}
This section introduces the main issue when proving Proposition \ref{EE} in a simple case. It allows to understand our strategy and the role of some of the assumptions in \textbf{(H)}. Furthermore, this presents in a simple form the effects of the degeneracy. By ``simple'', we mean that the assumptions  \textbf{(HL)} below hold true.\\

\textbf{Assumptions} \textbf{(HL)}. We say that hypotheses \textbf{(HL)} hold if \textbf{(H)} and \textbf{(HR)} hold with :
$F_1 \equiv 0_{\R^d}$, $\sigma \equiv {\rm Id}$, for all $(t,x_1,x_2)$ in $[0,T]\times \R^{d} \times \R^{d},\ F_2(t,x_1,x_2)=\bar{F}_2(x_2) + \Gamma_t x_1$. This implies in particular that for all $t$ in $[0,T]$, $\Gamma_t$ belongs to the convex subset $\mathcal{E}$ of $\rm{GL}_d(\R)$.\\

Under \textbf{(HL)}, the SDE (\ref{systemEDO}) becomes:
\begin{equation}\label{kolmog}
\left\lbrace\begin{array}{llll}
\d X^{1,t,x}_s= \d W_s,\qquad & X^{1,t,x}_t=x_1,\\
\d X^{2,t,x}_s = (\bar{F}_2(X^{2,t,x}_s) +\Gamma_s X^{1,t,x}_s) \d s, & X^{2,t,x}_t=x_2,\\
\end{array} \right.
\end{equation}
for all $t<s$ in $[0,T]^2$, $x$ in $\R^{2d}$ and admits a unique strong solution $X$.  We recall that the associated PDE is
\begin{equation}\label{PDELBH}
\left\lbrace\begin{array}{ll}
\partial_t u_i(t,x) + \mL u_i(t,x) = \phi_i(t,x), \quad \text{for }(t,x) \in [0,T] \times \mathbb{R}^{2d}\\
u_i(T,x)=0,\quad i=1,2,\\
\end{array}
\right.
\end{equation}
where $\mL$ is the generator of \eqref{kolmog}.\\

Our strategy to study the solution of the PDE \eqref{PDELBH} rests upon parametrix. This method is based on the following observation: in small time, the generator of the solution of an SDE with smooth and variable coefficients and the generator of the solution of the same SDE with fixed (frozen at some point) coefficients are ``closed’’. The variable generator is then seen as a perturbation of the frozen generator, which usually enjoys well known properties. 

Here, we know the explicit form of the fundamental solution of the frozen generator (which is the transition density of the solution of the frozen SDE).  Especially, we can prove that this fundamental solution and its derivatives admit Gaussian type bounds. Hence, thanks to the parametrix, we write the solution of the PDE \eqref{PDELBH} as a time-space integral of some perturbed kernel against this fundamental solution and we can study it.

We emphasize that the choice of the freezing point for the coefficients plays a central role: the perturbation done in the parametrix has to be of the order of the typical trajectories of the process associated to the frozen operator.

\subsection{The frozen system}\label{SZ1}

\textbf{Kolmogorov's example.} To understand how the frozen system could be, we go back to the work of Kolmogorov \cite{kolmogorov_zufallige_1934}  where the author studied the prototype system \eqref{systemEDO}. When $d=1$, Kolmogorov showed that the solution of $\d Y^1_s=\d W_s,\ \d Y^2_s = \alpha Y^1_s \d s$,  ($\alpha \neq 0$), with initial condition $(x_1,x_2)$ in $\R^2$, admits a density. Notably, this density is Gaussian and given by, for all $s$ in $(0,T]$ and $(y_1,y_2)$ in $\R^{2}$
\begin{equation}
p(0,x_1,x_2;s,y_1,y_2) = \frac{\sqrt{3}}{\alpha \pi s^2} \exp\left(-\frac{1}{2}\left|K_s^{-1/2}(y_1-x_1, y_2-x_2-s\alpha x_1)^* \right|^2\right),
\end{equation}	
with the following covariance matrix
\begin{equation}\label{COVMATKOKO}
K_s:=\begin{pmatrix}
s & (1/2)\alpha s^2 \\
(1/2)\alpha s^2 & (1/3)\alpha^2s^3 \\ 
\end{pmatrix}.
\end{equation} 

This example illustrates the behaviour of the system in small time: it is not diffusive. The first coordinate oscillates with fluctuations of order $1/2$, while the second one oscillates with fluctuations of order $3/2$. As a direct consequence, the transport of the initial condition of the first coordinate has a key role in the second one. This observation is crucial when freezing the coefficients.\\

\textbf{The frozen system.} As Kolmogorov's example suggests, we have to keep track of the transport of the initial condition when we freeze the coefficients. This allows us to reproduce a perturbation of the order of the typical trajectories of the frozen process. Then, we freeze the system \eqref{kolmog} along the forward transport $\theta_{\tau,s} = \left(\theta_{\tau,s}^1,\theta_{\tau,s}^2 \right)^*$, $s$ in $[\tau,T]$ that solves the ODE: 
\begin{equation}\label{c1:meanode}
\frac{\d}{\d s} \theta_{\tau,s} = \left(0_{\R^d},\ \bar{F}_2(\theta^2_{\tau,s}(\xi)) + \Gamma_s \theta^1_{\tau,s}(\xi)\right)^*,\ \theta_{\tau,\tau}(\xi)=\xi,
\end{equation}
for a given $(\tau,\xi)$ in $[0,T]\times \R^{2d}$ (we emphasize that in the regularized setting this ODE is well-posed). 
Hence, we obtain the following frozen system
\begin{equation}\label{kolmogf}
\left\lbrace\begin{array}{llll}
\d\bar{X}^{1,t,x}_s= \d W_s, \qquad & \bar{X}^{1,t,x}_t=x_1,\\
\d\bar{X}^{2,t,x}_s = \left(\bar{F}_2(\theta_{\tau,s}^2(\xi)) +  \Gamma_s \bar{X}^{1,t,x}_s\right) \d s, & \bar{X}^{2,t,x}_t=x_2,\\
\end{array} \right.
\end{equation}
for all $s$ in $(t,T]$. This is our candidate to approximate (\ref{kolmog})\footnote{Note that if $\tau>s$, we can extend the definition of $\theta$ and suppose that $ \forall (v>r)\in [0,T]^2, \forall \xi \in \R^{2d},\quad \theta_{v,r}(\xi)= 0$.}. Obviously, in order to reproduce the typical trajectories of the frozen process, the couple of variables $(\tau,\xi)$ in \eqref{c1:meanode} will be chosen as the initial data $(t,x)$ of the solution of the SDE \eqref{kolmogf}.

\subsection{Existence and Gaussian bound of the density of the frozen system}\label{SZ2}
In this case, the crucial point is the specific form of the covariance matrix $\bar{\Sigma}_{t,\cdot}$ of $\bar{X}^{t,x}_\cdot$ defined in \eqref{kolmogf}. For any  $s$ in $(t,T]$, standard computations show that
\begin{equation}\label{covar}
\bar{\Sigma}_{t,s}=\begin{pmatrix}
(s-t) & \int_t^s \int_t^r \Gamma_u \d u \d r\\
 \int_t^s  \int_t^r \Gamma_u \d u\d r &  \int_t^s  \left(\int_t^r \Gamma_u \d u\right)\left(\int_t^r \Gamma_u \d u\right)^* \d r\\
\end{pmatrix}.
\end{equation}

Therefore, the existence and the Gaussian estimates of the transition density of $\bar{X}^{t,x}_s$ stem from the control of the spectrum of $\bar{\Sigma}_{t,s}$. Such an investigation has been already done by  Delarue and Menozzi in \cite{delarue_density_2010}. The two following Lemmas shortly describe some of their results that are useful for us. The proofs are not given. For further details, we refer to Section 3 pp 18-24 of their paper. They prove that


\begin{lemme}\label{c1:lelemme}
Suppose that assumptions \textbf{(HL)} hold, then, a sufficient condition for the non-degeneracy of the variance matrix $\bar{\Sigma}_{t,s}$, $s$ in $(t,T]$ is given by
$$\det[\Gamma_r ]>0 \text{ for a.e. } r\in [t,s].$$
In that case, the solution of (\ref{kolmogf}) admits a transition density $\bar{q}$ given by, for all $s$ in $(t,T]$ and $y_1$, $y_2$ in $\R^d$:
\begin{equation}\label{TDFBH}
\bar{q}(t,x_1,x_2;s,y_1,y_2)= \frac{1}{(2\pi)^{d/2}} \det(\bar{\Sigma}_{t,s})^{-1/2} \exp \left(-\frac{1}{2}|\bar{\Sigma}^{-1/2}_{t,s}(y_1-x_1, y_2-m_{t,s}^{2,\tau,\xi}(x))^*  |^2\right),
\end{equation}
where 
$$m_{t,s}^{2,\tau,\xi}(x)=x_2 + \int_t^s\Gamma_r x_1 \d r  + \int_t^s \bar{F}_2(\theta_{\tau,r}^2(\xi))\d r,$$ 
and where $\bar{\Sigma}_{t,s}$ is the uniformly non-degenerate matrix given by \eqref{covar}.\\
\end{lemme}

From this expression, we can give the following Gaussian type estimate on the transition density of the solution of the SDE \eqref{kolmogf} and on its derivatives (the Gaussian bounds on the derivatives are not proven in the aforementioned work, but are proven in a more general case in the proof of Proposition \ref{estfundsol} below):
\begin{lemme}\label{lemmeLBH}
Suppose that assumptions \textbf{(HL)} hold, then, the transition kernel $\bar{q}$ and its derivatives admit Gaussian-type bounds: there exists a positive constant $c$ depending only on known parameters in \textbf{(H)} such that for all $\xi$ in $\R^{2d}$:
\begin{eqnarray}\label{oto}
&&\left| D_{x_1}^{N^{x_1}} D_{x_2}^{N^{x_2}} D_{y_1}^{N^{y_1}} \bar{q} (t,x_1,x_2;s,y_1,y_2) \right| \\
&&\leq (s-t)^{-[3N^{x_2} + N^{x_1} + N^{y_1}]/2} \frac{c}{(s-t)^{2d}} \exp \left( -c \left(\frac{\big|y_1-x_1\big|^2}{s-t}+ \frac{\big|y_2-m^{2,\tau,\xi}_{t,s}(x)\big|^2}{(s-t)^{3}} \right)\right),\nonumber
\end{eqnarray}
for all $s$ in $(t,T]$, $y_1$, $y_2$ in $\R^d$ and any $N^{x_1},N^{x_2},N^{y_1}$ less than 2. 
\end{lemme}

We emphasize that the constant $c$ that appears in the exponential in estimate \eqref{oto} does not depend on $\Gamma$, as suggested by Lemma \ref{c1:lelemme}. This uniform control is not obvious and is related to the ``closed convex'' assumption \textbf{(H3-b)}.

If the control is not uniform, Delarue and Menozzi show in \cite{delarue_density_2010} (see Example 3.5) that one can find a sequence of matrix $(\Gamma_1^n)_{n\geq 0}$ with positive constant determinant such that $\det[\Sigma^n_{0,1}]$ converges towards 0 as $n$ tends to the infinity. The crucial point in their example is that the sequence of functions $(t\in [0,1]\mapsto \Gamma_t^n)_{n\geq 0}$ weakly converges towards 0. Hence, to overcome this problem, the authors need some closure for the weak topology. The closed convex assumption \textbf{(H3-b)} allows them to obtain compactness for the weak topology.\\


Note that for all $s$ in $[t,T]$, the mean $(x_1,m^{2,t,x}_{t,s}(x))$ of $\bar{X}_s^{t,x}$ satisfies the ODE \eqref{c1:meanode} with initial data $(\tau,\xi)=(t,x)$. Since this equation admits a unique solution under \textbf{(HL)}, we deduce that for all $s$ in $[t,T]$, the forward transport function defined by \eqref{c1:meanode} with the starting point $x$ as initial condition is equal to the mean: $\theta_{t,s}(x)=m^{t,x}_{t,s}(x)$.

Finally, as we will show in the proof of Proposition \ref{estfundsol} in Section \ref{SOTLSOP}, the transition density $\bar{q}$ enjoys the following symmetry  property:
\begin{equation}\label{c1:centeringZA}
\forall (t<s,x,y)\in [0,T]^2\times \R^{d} \times \R^{d}, D_{x_2}\bar{q}(t,x;s,y) = -D_{y_2}\bar{q}(t,x;s,y).
\end{equation} 
This plays a crucial role in the proof of the Lipschitz estimates of the solution $u$ of the PDE \eqref{PDELBH} and of its derivative $D_{1}u$ .

\subsection{Representation of the solution by parametrix}\label{SZ3}
Let $\bar{\mL}^{\tau,\xi} := (1/2)\Delta_{x_1} + \left[\bar{F}_2(\theta_{\tau,t}^2(\xi)) + \Gamma_t x_1\right] \cdot D_{x_2}$ be the generator of the frozen process $\bar{X}$. We can write the PDE \eqref{PDELBH} as
\begin{equation*}\left\lbrace\begin{array}{ll}
\partial_t u_i(t,x) + \bar{\mL}^{\tau,\xi} u_i(t,x) = \phi_i(t,x) + (\bar{\mL}^{\tau,\xi}-\mL)u_i(t,x) , \ \text{for }(t,x) \in [0,T) \times \mathbb{R}^{2d}\\ u_i(T,x)=0,\quad i=1,2.
\end{array}
\right.
\end{equation*} 
Since $\bar{q}$ is the fundamental solution of $\bar{\mL}^{\tau,\xi}$ we have, for all $(t,x)$ in $[0,T]\times \R^{2d}$, that the unique solution of this PDE reads:
\begin{eqnarray*}
u_i(t,x) &=& \int_t^T \int_{\R^{2d}} \bigg\{\phi_i(s,y) - [\bar{F}_2(y_2) -\bar{F}_2(\theta_{\tau,s}^2(\xi)) ] \cdot D_{y_2}u_i(s,y)\bigg\}\bar{q}(t,x;s,y) \d y\d s,
\end{eqnarray*}
for $i=1,2$ and for all $(\tau,\xi)$ in $[0,T]\times\R^{2d}$. We emphasize that $u$ is the unique solution of \eqref{PDELBH} and does not depend on the choice of the freezing point ``$(\tau,\xi)$''. For every $(t,x)$ in $[0,T]\times \R^{2d}$, we can then choose $\tau$ as the current evaluation time and then write
\begin{eqnarray}\label{KolPDEsol}
u_i(t,x) &=& \int_t^T \int_{\R^{2d}} \bigg\{\phi_i(s,y) - [\bar{F}_2(y_2) -\bar{F}_2(\theta_{t,s}^2(\xi)) ] \cdot D_{y_2}u_i(s,y)\bigg\}\bar{q}(t,x;s,y) \d y\d s.
\end{eqnarray}
From now on, for every given time $t$ in $[0,T]$, we let $\tau=t$.
\subsection{A priori estimates}\label{SZ4}
In order to prove the bounds of Proposition \ref{EE}, we need to obtain estimates of the supremum norm of the first and second order derivatives of the $u_i$, $i=1,2$. Having in mind to invert the differentiation and integral operators, we see any differentiation of $u_i$ as an integral of a certain function against the derivative of the degenerate Gaussian kernel $\bar{q}$. 

As shown in Lemma \ref{lemmeLBH}, the differentiation of this kernel generates a time-singularity. Each differentiation of the transition kernel w.r.t. the first component gives a time-singularity of order $1/2$ while the differentiation w.r.t. the second component gives a time-singularity of order $3/2$. 

The main issue consists in smoothing this singularity by using the regularity of the coefficients assumed in \textbf{(H)} together with Gaussian decay in $\bar{q}$ by letting the freezing point $\xi$ be the starting point of the process.

Let us illustrate the computations with the worst case in Proposition \ref{EE}, that is, the cross derivative  $D_{1}D_{2}u_i$ (which gives a time singularity of order 2, see Lemma \ref{lemmeLBH}). In order to invert the integral and the differentiation operators we have to show that the derivative of the time integrand in \eqref{KolPDEsol} is suitably bounded. When it is evaluated at point $\xi=x$ we show in the sequel that this is indeed the case. Let $(t,x)$ in $[0,T]\times \R^{2d}$, for $i$ in $\{1,2\}$, we denote by $\mathcal{I}_i(s,x)$, $s$ in $(t,T]$, the  time integrand in \eqref{KolPDEsol}. By switching the differentiation and the (space) integral operator we have:
\begin{eqnarray}
 &&D_{x_{1}} D_{x_{2}} \mI_i(s,x)\label{diffeu}\\
 && \quad =  \int_{\R^{2d}} \bigg\{\phi_i(s,y) - [\bar{F}_2(y_2) -\bar{F}_2(\theta_{t,s}^2(\xi)) ] \cdot D_{y_2}u_i(s,y)\bigg\}  D_{x_{1}} D_{x_{2}} \bar{q}(t,x;s,y) \d y \notag\
\end{eqnarray}

We deal with the right hand side of \eqref{diffeu}. Let us first deal with the first term in the time-integrand. For all $s$ in $(t,T]$ we have:
\begin{eqnarray*}
&&\int_{\R^{2d}} \phi_i(s,y_1,y_2) D_{x_{1}} D_{x_{2}} \bar{q}(t,x_1,x_2;s,y_1,y_2) \d y_1\d y_2\\
&&=\int_{\R^{2d}} \left(\phi_i(s,y_1,y_2)-\phi_i(s,y_1,\theta_{t,s}^2(\xi)) \right) D_{x_{1}} D_{x_{2}} \bar{q}(t,x_1,x_2;s,y_1,y_2) \d y_1\d y_2\\
&&\quad +\int_{\R^{2d}} \phi_i(s,y_1,\theta_{t,s}^2(\xi)) D_{x_{1}} D_{x_{2}} \bar{q}(t,x_1,x_2;s,y_1,y_2)\d y_1\d y_2,
 \end{eqnarray*}
where, thanks to the symmetry \eqref{c1:centeringZA} and an integration by parts argument, the last term in the right hand side is equal to 0. In the sequel, we refer to this argument as \emph{the centering argument}. Combining this argument and the estimate of $D_{x_{1}} D_{x_{2}}\bar{q}$ in Lemma \ref{lemmeLBH}, we have, for all $s$ in $(t,T]$:
\begin{eqnarray*}
&&\left|\int_{\R^{2d}} \left(\phi_i(s,y_1,y_2)-\phi_i(s,y_1,\theta_{t,s}^2(\xi)) \right) D_{x_{1}}D_{x_{2}} \bar{q}(t,x_1,x_2;s,y_1,y_2) \d y_1\d y_2\right|\\ 
&& \leq \int_{\R^{2d}} \Bigg\{ (s-t)^{-2}\left|\phi_i(s,y_1,y_2)-\phi_i(s,y_1,\theta_{t,s}^2(\xi)) \right| \\
&& \qquad \times \frac{c}{(s-t)^{2d}} \exp \left( -c \left(\frac{\big|y_1-x_1\big|^2}{s-t}+ \frac{\big|y_2-m^{2,t,\xi}_{t,s}(x)\big|^2}{(s-t)^{3}} \right)\right) \Bigg\} \d y_1\d y_2,
\end{eqnarray*}
where $c$ depends only on known parameters in \textbf{(H)}. By using the H\"{o}lder regularity of $\phi_i$ assumed in \textbf{(H1)} we have
\begin{eqnarray}
&&\left|\int_{\R^{2d}}(s-t)^{-2} \left|\phi_i(s,y_1,y_2) - \phi_i(s,y_1,\theta_{t,s}^2(\xi))\right|\bar{q}(t,x_1,x_2;s,y_1,y_2)\d y_1\d y_2 \right| \notag\\
&& \quad \leq C \int_{\R^{2d}} \Bigg\{(s-t)^{-2} |y_2-\theta_{t,s}^2(\xi)|^{\beta_i^2} \notag\\
&& \qquad \times \frac{c}{(s-t)^{2d}} \exp \left( -c \left(\frac{\big|y_1-x_1\big|^2}{s-t}+ \frac{\big|y_2-m^{2,t,\xi}_{t,s}(x)\big|^2}{(s-t)^{3}} \right)\right)\Bigg\}\d y_1\d y_2.\label{lexpoc}
\end{eqnarray}  
Now, we use the Gaussian off-diagonal decay of $\bar{q}$ to smooth the time-singularity: by letting $\xi=x$ (and then $\theta_{t,s}^2(x)=m^{2,t,x}_{t,s}(x)$), for all positive $\eta$, there exists a positive constant $\bar{C}$ such that\footnote{By using the inequality: $\forall \eta>0,\ \forall q>0,\ \exists \bar{C}>0 \text{ s.t. } \forall \sigma >0,\ \sigma^{q} e^{-\eta\sigma}\leq \bar{C} $.} 
$$\left(\frac{|y_2-m^{2,t,x}_{t,s}(x)|}{(s-t)^{3/2}}\right)^{\beta_i^2} \times \exp \left( -\eta\left(\frac{\big|y_1-x_1\big|^2}{s-t}+ \frac{\big|y_2-m^{2,t,x}_{t,s}(x)\big|^2}{(s-t)^{3}} \right)\right)\leq \bar{C},$$
where $\bar{C}$ depends on $\eta$ and $\beta_i^2$ only. Thus, by damaging the constant $c$ in the exponential in \eqref{lexpoc}, we obtain the following estimate:
\begin{eqnarray}\label{diagdec}
&& \left|\int_{\R^{2d}} (s-t)^{-2} \left|\phi_i(s,y_1,y_2) - \phi_i(s,y_1,\theta_{t,s}^2(x))\right|\bar{q}(t,x_1,x_2;s,y_1,y_2)\d y_1\d y_2 \right| \\
&& \quad \leq C' \int_{\R^{2d}} (s-t)^{-2+3\beta_i^2/2} \frac{c}{(s-t)^{2d}} \exp \left( -c' \left(\frac{\big|y_1-x_1\big|^2}{s-t}+ \frac{\big|y_2-m^{2,t,x}_{t,s}(x)\big|^2}{(s-t)^{3}} \right)\right)\d y_1\d y_2\notag.
\end{eqnarray} 

Therefore, by choosing the value of $\beta_i^2$ strictly greater than $2/3$, the singularity  $(s-t)^{-2+3\beta_i^2/2}$ becomes integrable. By applying the same procedure (without centering) with the term
$$ \int_{\R^{2d}} \bigg\{[\bar{F}_2(y_2) -\bar{F}_2(\theta_{t,s}^2(\xi)) ] \cdot D_{y_2}u_i(s,y)\bigg\}\bar{q}(t,x;s,y)\d y ,$$
we can deduce that
\begin{equation}\label{LBH:crossesti}
\left\|D_{{1}}D_{{2}}u_i\right\|_{\infty} \leq C''(T^{-1+3\beta_i^2/2}+T^{-1+3\beta_2^2/2}\left\|D_{{2}}u_i\right\|_{\infty}).
\end{equation}

The main problem here is that the supremum norm of $D_{{2}}u_i$ appears in the bound, so that we also have to estimate this quantity in term on known parameters in \textbf{(H)} to overcome the problem. It is well seen that the same arguments lead to
$$\left\|D_{{2}}u_i\right\|_{\infty} \leq C'''(T^{(-1+3\beta_i^2)/2}+T^{(-1+3\beta_2^2)/2}\left\|D_{{2}}u_i\right\|_{\infty}).$$
By choosing $T$ sufficiently small (\emph{e.g.} such that $C'''T^{(-1+3\beta_2^2)/2}$ is less than $1/2$) we obtain that
\begin{equation}\label{circular}
\left\|D_{{2}}u_i\right\|_{\infty} \leq 2C'''T^{(-1+3\beta_i^2)/2}.
\end{equation}
We refer to this argument as the \emph{circular argument} in the following.  By plugging this bound in \eqref{LBH:crossesti} and by applying the same strategy with $D_{1}u_i$ and $D^2_{1}u_i$, Proposition \ref{EE} under \textbf{(HL)} follows for $T$ less than $\mathcal{T} := (1/(2C'''))^{2/(-1+3\beta_2^2)}$.\\
 
From this discussion, one can also see the specific choice of the freezing curve as the one that matches the off-diagonal decay of the exponential in $\bar{q}$ when $\xi=x$.

\section{Mathematical tools}\label{SOTLSOP}
In this section, we introduce the ingredients for the proof of Proposition \ref{EE}.

\subsection{The frozen system}\label{LSAP}
Consider the frozen system:
\begin{equation}\label{LS}
\begin{array}{ll}
 \d\tilde{X}^{1,t,x}_s = F_1(s,\theta_{\tau,s}(\xi)) \d s + \sigma(s,\theta_{\tau,s}(\xi)) \d W_s\\
 \d\tilde{X}^{2,t,x}_s = \left[F_2(s,\theta_{\tau,s}(\xi)) + D_{1}F_2(s,\theta_{\tau,s}(\xi))(\tilde{X}^{1,t,x}_s-\theta^1_{\tau,s}(\xi)) \right]\d s
\end{array}
\end{equation}
for all $s$ in $(t,T]$, any $t$ in $[0,T]$, and for any initial condition $x$ in $\R^{2d}$ at time $t$ and any $(\tau,\xi) \in [0,T]\times \R^{2d}$, linearized around the transport $(\theta_{\tau,s}(\xi))_{\tau \leq s \leq T}$ defined\footnote{Note again that if $\tau>s$, we can extend the definition of $\theta$ and suppose that $ \forall (v>r)\in [0,T]^2, \forall \xi \in \R^{2d},\quad \theta_{v,r}(\xi)= 0$.} by
\begin{equation}\label{thetadef}
\frac{\d}{\d s}\theta_{\tau,s}(\xi) = F(s,\theta_{\tau,s}(\xi)),\quad \theta_{\tau,\tau}(\xi)=\xi.
\end{equation}
The following Proposition holds:

\begin{proposition}\label{estfundsol}
Suppose that assumptions \textbf{(HR)} hold, then:

(i) There exists a unique (strong) solution of \eqref{LS}  with mean 
$$(m^{\tau,\xi}_{t,s})_{t \leq s \leq T} =(m^{1,\tau,\xi}_{t,s},m^{2,\tau,\xi}_{t,s})_{t \leq s \leq T},  $$
where
\begin{eqnarray}\label{meanGauss}
&&m^{1,\tau,\xi}_{t,s}(x) =  x_1 + \int_t^s F_1(r,\theta_{\tau,r}(\xi)) \d r,\\
&&m^{2,\tau,\xi}_{t,s}(x) = x_2+\int_t^s  \bigg[F_2(r,\theta_{\tau,r}(\xi)) + D_{1}F_2(r,\theta_{\tau,r}(\xi))(x_1-\theta_{\tau,r}^1(\xi)) \nonumber\\
&& \hphantom{m^{2,\xi}_{\tau,s}(x)} \quad +  D_{1}F_2(r,\theta_{\tau,r}(\xi))\int_t^r F_1(v,\theta_{\tau,v}(\xi))\d v \bigg] \d r,\nonumber
\end{eqnarray} 
and uniformly non-degenerate covariance matrix $(\tilde{\Sigma}_{t,s})_{t\leq s \leq T}$:

\begin{equation}\label{covmatrice}
\tilde{\Sigma}_{t,s}=\begin{pmatrix}
\int_t^s \sigma \sigma^*(r,\theta_{\tau,r}(\xi))\d r & \int_t^s R_{r,s}(\tau,\xi) \sigma \sigma^*(r,\theta_{\tau,r}(\xi))\d r\\
 \int_t^s \sigma \sigma^*(r,\theta_{\tau,r}(\xi))R^*_{r,s}(\tau,\xi) \d r &  \int_t^s R_{t,r}(\tau,\xi)\sigma \sigma^*(r,\theta_{\tau,r}(\xi)) R_{t,r}^*(\tau,\xi) \d r\\
\end{pmatrix},					
\end{equation}
where:
\begin{equation*}
R_{t,r}(\tau,\xi)=\left[\int_t^r D_{1}F_2(v,\theta_{\tau,v}(\xi))\d v\right],\quad t\leq r\leq s\leq T.
\end{equation*}

(ii) This solution is a Gaussian process with transition density:
\begin{eqnarray}\label{gtd}
\tilde{q}(t,x_1,x_2;s,y_1,y_2)  = \frac{3^{d/2}}{(2\pi)^{d/2}}  (\det[\tilde{\Sigma}_{t,s}])^{-1/2} \exp \left( -\frac{1}{2}|\tilde{\Sigma}_{t,s}^{-1/2} (y_1-m^{1,\tau,\xi}_{t,s}(x), y_2-m^{2,\tau,\xi}_{t,s}(x))^*  |^2\right),
\end{eqnarray}
for all $s$ in $(t,T]$.

(iii) This transition density $\tilde{q}$ is the fundamental solution of the PDE driven by $\tilde{\mL}^{\tau,\xi}$ and given by:
\begin{eqnarray}
\tilde{\mL}^{\tau,\xi} &:=& \frac{1}{2} Tr\left[a(t,\theta_{\tau,t}(\xi))D^2_{x_1}\right] +  \left[ F_1(t,\theta_{\tau,t}(\xi))\right] \cdot D_{x_1}  \nonumber\\
&& \quad + \left[F_2(t,\theta_{\tau,t}(\xi))+ D_{1} F_2(t,\theta_{\tau,t}(\xi))\left(x_1 - \theta^1_{\tau,t}(\xi)\right)\right] \cdot D_{x_2}. \label{frozgen}
\end{eqnarray}

(iv) There exist two positive constants $c$ and $C$, depending only on known parameters in \textbf{(H)}, such that
\begin{equation}\label{c1defqc}
\tilde{q}(t,x_1,x_2;s,y_1,y_2) \leq C\hat{q}_{c}(t,x_1,x_2;s,y_1,y_2),
\end{equation}
where
\begin{eqnarray*}
\hat{q}_{c}(t,x_1,x_2;s,y_1,y_2)= \frac{c}{(s-t)^{2d}}\exp \left( -c \left(\frac{\big|y_1-m^{1,\tau,\xi}_{t,s}(x)\big|^2}{s-t}+ \frac{\big|y_2-m^{2,\tau,\xi}_{t,s}(x)\big|^2}{(s-t)^{3}} \right)\right),
\end{eqnarray*}
and
\begin{eqnarray}\label{estidertranker}
\left| D_{x_1}^{N^{x_1}} D_{x_2}^{N^{x_2}} D_{y_1}^{N^{y_1}} \tilde{q}(t,x_1,x_2;s,y_1,y_2)\right|  \leq C (s-t)^{-[3N^{x_2} + N^{x_1} + N^{y_1}]/2} \hat{q}_c(t,x_1,x_2;s,y_1,y_2),
\end{eqnarray}
for all $s$ in $(t,T]$ and any integers $N^{x_1},N^{x_2},N^{y_1}$ less than 2.\\
\end{proposition}

\begin{proof}
(i) First of all, note that, under \textbf{(HR)}, the ODE: $[\d/\d s]\theta_{\tau,s}(\xi) = F(s,\theta_{\tau,s}(\xi))$, $\theta_{\tau,\tau}(\xi)=\xi$ admits a unique solution and that \eqref{LS} admits a unique strong solution $\tilde{X}$. By rewriting \eqref{LS} in integral form and by plugging the obtained representation of $\tX^1$ in $\tX^2$, it is easily seen that the expressions of the mean \eqref{meanGauss} and the variance \eqref{covmatrice} follow from the stochastic Fubini Theorem and standard computations. The uniform non-degeneracy of $(\tilde{\Sigma}_{t,s})_{t < s \leq T}$ arises from assumptions \textbf{(H)} and Proposition 3.1 in \cite{delarue_density_2010}.\\

(ii)-(iii) These assertions result from standard computations.\\

(iv) For all $s$ in $(t,T]$, we know from Proposition 3.1 in \cite{delarue_density_2010} that the matrix $\tilde{\Sigma}_{t,s}$ is symmetric and uniformly non-degenerate. Besides, from Proposition 3.4 in \cite{delarue_density_2010} there exists a constant $C$ depending only on known parameters in \textbf{(H)} such that for all $s$ in $(t,T]$, for all $(x,y,\xi)$ in $\R^{2d} \times \R^{2d} \times \R^{2d} $,
\begin{eqnarray*}
&&-\left[\tilde{\Sigma}_{t,s}^{-1} (y_1-m^{1,\tau,\xi}_{t,s}(x), y_2-m^{2,\tau,\xi}_{t,s}(x))^*\right] \cdot \left[(y_1-m^{1,\tau,\xi}_{t,s}(x), y_2-m^{2,\xi}_{t,s}(x))^*\right] \\
&& \qquad \leq - C \left[ \left(\frac{y_1-m^{1,\tau,\xi}_{t,s}(x)}{(s-t)^{1/2}}, \frac{y_2-m^{2,\tau,\xi}_{t,s}(x)}{(s-t)^{3/2}}\right)^* \right] \cdot \left[\left(\frac{y_1-m^{1,\tau,\xi}_{t,s}(x)}{(s-t)^{1/2}}, \frac{y_2-m^{2,\tau,\xi}_{t,s}(x)}{(s-t)^{3/2}}\right)^* \right].
\end{eqnarray*}

For $i,j=1,2$, let $[\tilde{\Sigma}^{-1}_{t,s}]_{i,j}$ denote the block of size $d\times d$ of the matrix $\tilde{\Sigma}^{-1}_{t,s}$ at the $(i-1)d+1, (j-1)d+1$ rank. We can deduce from (\ref{covmatrice}) that there exists a positive constant $C$ depending only on known parameters in \textbf{(H)} such that (we also refer the reader to Lemma 3.6 and to the proof of Lemma 5.5 in \cite{delarue_density_2010} for more details), for all $s$ in $(t,T]$, for all $\zeta$ in $\R^d$:
\begin{equation}\label{rereg}
\begin{array}{llll}
&\left|[\tilde{\Sigma}_{t,s}^{-1}]_{1,1}\zeta\right| \leq C (s-t)^{-1} \left|\zeta\right|,\\
&\left|[\tilde{\Sigma}_{t,s}^{-1}]_{1,2}\zeta\right| + \left|[\tilde{\Sigma}_{t,s}^{-1}]_{2,1}\zeta\right|\leq C (s-t)^{-2} \left|\zeta\right|,\\
&\left|[\tilde{\Sigma}_{t,s}^{-1}]_{2,2}\zeta\right|  \leq C (s-t)^{-3} \left|\zeta\right|,
\end{array}
\end{equation}
hence, $\tilde{\Sigma}^{-1}_{t,.}$ has the same structure as $K^{-1}_{.-t}$ in \eqref{COVMATKOKO}.\\

Now, we compute the derivatives w.r.t. each component and estimate it with the help of \eqref{rereg}. Let $(t<s,x,y)$ in $[0,T]^2 \times \R^{2d} \times \R^{2d}$, we have:
\begin{eqnarray*}
&&|D_{x_2} \tilde{q}(t,x_1,x_2;s,y_1,y_2)|\\
&&=  \left| \left( -2 [\tilde{\Sigma}^{-1}_{t,s}]_{2,1}(y_1-m^{1,\tau,\xi}_{t,s}(x)) -2 [\tilde{\Sigma}^{-1}_{t,s}]_{2,2}(y_2-m^{2,\tau,\xi}_{t,s}(x)) \right)\tilde{q}(t,x_1,x_2;s,y_1,y_2)\right|\\ 
&&\leq C (s-t)^{-3/2}\left(\left| \frac{(y_1-m^{2,\tau,\xi}_{t,s}(x))}{(s-t)^{1/2}}\right| + \left|\frac{(y_2-m^{2,\tau,\xi}_{t,s}(x))}{(s-t)^{3/2}} \right| \right) \tilde{q}(t,x_1,x_2;s,y_1,y_2)\\ 
&&\leq C' (s-t)^{-3/2} \hat{q}_c(t,x_1,x_2;s,y_1,y_2).
\end{eqnarray*}
Note that the symmetry $D_{x_2} \tilde{q}(t,x_1,x_2;s,y_1,y_2)=-D_{y_2} \tilde{q}(t,x_1,x_2;s,y_1,y_2)$ holds. Now, we have

\begin{eqnarray*}
&&|D_{y_1} \tilde{q}(t,x_1,x_2;s,y_1,y_2)| \\
&&= \left| \left(  2[\tilde{\Sigma}^{-1}_{t,s}]_{1,1}(y_1-m^{1,\tau,\xi}_{t,s}(x))  + 2 [\tilde{\Sigma}^{-1}_{t,s}]_{1,2}(y_2-m^{2,\tau,\xi}_{t,s}(x))\right)\tilde{q}(t,x_1,x_2;s,y_1,y_2) \right|\\
&&\leq C (s-t)^{-1/2} \hat{q}_c(t,x_1,x_2;s,y_1,y_2).\\
\end{eqnarray*}
Unfortunately, the transport of the initial condition of the diffusive component in the degenerate component breaks the symmetry and $D_{x_1} \tilde{q}(t,x_1,x_2;s,y_1,y_2)\neq -D_{y_1} \tilde{q}(t,x_1,x_2;s,y_1,y_2)$. Indeed

\begin{eqnarray*}
D_{x_1} \tilde{q}(t,x_1,x_2;s,y_1,y_2) &=& \bigg( -2[\tilde{\Sigma}^{-1}_{t,s}]_{1,1}(y_1-m^{1,\tau,\xi}_{t,s}(x)) -2[\tilde{\Sigma}^{-1}_{t,s}]_{1,2}(y_2-m^{2,\tau,\xi}_{t,s}(x))\\
 && - 2[\tilde{\Sigma}^{-1}_{t,s}]_{1,2}\left[\left(R_{t,s}(\tau,\xi)\right)(y_1-m^{1,\tau,\xi}_{t,s}(x))\right]  \\
&&  - 2[\tilde{\Sigma}^{-1}_{t,s}]_{2,2} \left[\left(R_{t,s}(\tau,\xi)\right)(y_2-m^{2,\tau,\xi}_{t,s}(x))\right]\bigg)\tilde{q}(t,x_1,x_2;s,y_1,y_2).
\end{eqnarray*}
Since the term  $R_{t,s}(\xi)$ is of order $(s-t)$ (this is the transport of the initial condition from time $t$ to $s$), we deduce that

\begin{eqnarray*}
|D_{x_1} \tilde{q}(t,x_1,x_2;s,y_1,y_2)| &\leq & C(s-t)^{-1/2}  \Bigg\{ \left|  \frac{(y_1-m^{2,\tau,\xi}_{t,s}(x))}{(s-t)^{1/2}}\right|  + \left|\frac{(y_2-m^{2,\tau,\xi}_{t,s}(x))}{(s-t)^{3/2}} \right|  \\
&&+\left| \frac{(y_1-m^{2,\tau,\xi}_{t,s}(x))}{(s-t)^{1/2}}\right|   + \left|\frac{(y_2-m^{2,\tau,\xi}_{t,s}(x))}{(s-t)^{3/2}} \right|  \Bigg\}  \tilde{q}(t,x_1,x_2;s,y_1,y_2)\\
&\leq & C' (s-t)^{-1/2} \hat{q}_c(t,x_1,x_2;s,y_1,y_2).
\end{eqnarray*}
Finally,

\begin{eqnarray*}
&&D^2_{x_1} \tilde{q}(t,x_1,x_2;s,y_1,y_2)\\
 &&=   \left( -2[\tilde{\Sigma}^{-1}_{t,s}]_{1,1}D_{x_1}m^{1,\tau,\xi}_{t,s}(x)) -2[\tilde{\Sigma}^{-1}_{t,s}]_{1,2}D_{x_1}m^{2,\tau,\xi}_{t,s}(x) - 2[\tilde{\Sigma}^{-1}_{t,s}]_{1,2}\left[\left(R_{t,s}(\tau,\xi)\right)D_{x_1}m^{1,\tau,\xi}_{t,s}(x)\right] \right.\\
&& \qquad \left.   - 2[\tilde{\Sigma}^{-1}_{t,s}]_{2,2} \left[\left(R_{t,s}(\tau,\xi)\right)D_{x_1}m^{2,\tau,\xi}_{t,s}(x)\right]\right) \tilde{q}(t,x_1,x_2;s,y_1,y_2)\\
&& \quad +\left( -2[\tilde{\Sigma}^{-1}_{t,s}]_{1,1}(y_1-m^{1,\tau,\xi}_{t,s}(x)) -2[\tilde{\Sigma}^{-1}_{t,s}]_{1,2}(y_2-m^{2,\tau,\xi}_{t,s}(x))  \right.\\
&& \qquad - 2[\tilde{\Sigma}^{-1}_{t,s}]_{1,2}\left[\left(R_{t,s}(\tau,\xi)\right)(y_1-m^{1,\tau,\xi}_{t,s}(x))\right]\\
&& \qquad \left.   - 2[\tilde{\Sigma}^{-1}_{t,s}]_{2,2} \left[\left(R_{t,s}(\tau,\xi)\right)(y_2-m^{2,\tau,\xi}_{t,s}(x))\right]\right)^2\tilde{q}(t,x_1,x_2;s,y_1,y_2).
\end{eqnarray*}
Note that, from (\ref{meanGauss}) we have $D_{x_1}m^{\tau,\xi}_{t,s}(x)=({\rm Id}, R_{t,s}(\tau,\xi))^*$, so that,

\begin{eqnarray*}
&&|D^2_{x_1} \tilde{q}(t,x_1,x_2;s,y_1,y_2)| \leq C (s-t)^{-1} \hat{q}_c(t,x_1,x_2;s,y_1,y_2).
\end{eqnarray*}

The other derivatives can be deduced from these computations and estimate (\ref{estidertranker}) follows.
\end{proof}

\begin{remarque}
From this proof, one can deduce that the symmetry $D_{x_2} \tilde{q}=-D_{y_2} \tilde{q}$ holds. Therefore, by an integration by parts argument, for all $t$ in $[0,T]$, all $s$ in $[t,T]$ and  $y_1,x_1,x_2$ in $\R^d$,
\begin{eqnarray}
&&\int_{\R^d} D_{x_2}\tilde{q}(t,x_1,x_2;s,y_1,y_2) dy_2 =0.\label{argcentre}
\end{eqnarray}
This argument is very useful in the sequel.
\end{remarque}

\subsection{Definitions and rules of calculus}
We introduce some definitions and rules of computations that will be useful in the following section. Let us begin by the following definition:

\begin{definition}\label{chap1:perturbed} For all $\zeta$ in $\R^{2d}$ we denote by $\Delta(\zeta)$ the perturbation operator around $\zeta$ acting on any function $f$ from $[0,T]\times \R^{2d}$ as follows:
$$\forall (s,y) \in  [0,T]\times \R^{2d},\ \Delta(\zeta)f(s,y)=f(s,y)-f(s,\zeta),$$
and for $j=1,2$, we denote by $\Delta^j(\zeta)$ the perturbation operator around $\zeta_j$ acting on any function $f$ from $[0,T]\times \R^{2d}$ as follows:
$$\forall (s,y) \in  [0,T]\times \R^{2d},\ \Delta^1(\zeta) f(s,y_1,y_2)=f(s,y_1,\zeta_2)-f(s,\zeta_1,\zeta_2),$$
and
$$\forall (s,y) \in  [0,T]\times \R^{2d},\ \Delta^2(\zeta)f(s,y_1,y_2)=f(s,y_1,y_2)-f(s,y_1,\zeta_2).$$

Especially, the notation $\Delta^j(\zeta)y_j$ stands for $y_j-\zeta_j$.
\end{definition}

Given this definition we can give a generic \emph{centering argument}, as introduced in Subsection \ref{SZ4} in the linear and Brownian heuristic:

\begin{claim}{}\label{c1center}
Let $\tq$ be the function defined by \eqref{gtd} in Proposition \ref{estfundsol} and let $f$ and $g$ be two continuous functions defined on $[0,T]\times \R^{2d}$. For all $N^1$ $N^2$ in $\mathbb{N}$, for all $t< s$ in $[0,T]^2$, $x$ in $\R^{2d}$ and $\zeta$ in $\R^{2d}$ we have that:
\begin{enumerate}[(a)]
\item $\displaystyle D_{x_{1}}^{N^1}D_{x_{2}}^{N^2} \int_{\R^{2d}}  f(s,y) \tq(t,x;s,y)  \d y  =D_{x_{1}}^{N^1}D_{x_{2}}^{N^2}  \int_{\R^{2d}}  \Delta(\zeta) f(s,y) \tq(t,x;s,y)  \d y ,$\label{c1centeresti1}\\
if $N_1+N_2>0$,
\item $\displaystyle D_{x_{1}}^{N^1}D_{x_{2}}^{N^2}  \int_{\R^{2d}} f(s,y)\tq(t,x;s,y)  \d y = D_{x_{1}}^{N^1}D_{x_{2}}^{N^2}   \int_{\R^{2d}} \Delta^2(\zeta) f(s,y) \tq(t,x;s,y)  \d y ,$\label{c1centeresti2}\\
if $N_2>0$,
\item $\displaystyle D_{x_{1}}^{N^1}D_{x_{2}}^{N^2} \int_{\R^{2d}}  \Delta^1(\zeta) f(s,y) g(s,y)\tq(t,x;s,y)  \d y $\\
$\displaystyle\qquad  =D_{x_{1}}^{N^1}D_{x_{2}}^{N^2} \int_{\R^{2d}}  (\Delta^1(\zeta) f(s,y))( \Delta^2(\zeta) g(s,y)) \tq(t,x;s,y)  \d y $,\label{c1centeresti3} if $N_2>0$.
\end{enumerate}
\end{claim}

\begin{proof}
Let $f$ and $g$ be defined as in Claim \ref{c1center} and let $t< s$ in $[0,T]^2$, $x$ in $\R^{2d}$. We have, by Definition \ref{chap1:perturbed}:
\begin{eqnarray*}
D_{x_{1}}^{N^1}D_{x_{2}}^{N^2}  \int_{\R^{2d}} f(s,y)\tq(t,x;s,y)  \d y &=& D_{x_{1}}^{N^1}D_{x_{2}}^{N^2} \int_{\R^{2d}}  \Delta(\zeta) f(s,y) \tq(t,x;s,y)  \d y  \nonumber\\
&&+ D_{x_{1}}^{N^1}D_{x_{2}}^{N^2} \int_{\R^{2d}} f(s,\zeta) \tq(t,x;s,y)  \d y, \nonumber
\end{eqnarray*}
for all $\zeta$ in $\R^{2d}$. The last term in the right hand side is equal to 0 since it does not depend on $x$ after integrating. This concludes the proof of \eqref{c1centeresti1}. Now, we prove \eqref{c1centeresti3}. For all $\zeta$ in $\R^{2d}$ we have 
\begin{eqnarray*}
&&D_{x_{1}}^{N^1}D_{x_{2}}^{N^2} \int_{\R^{2d}} (\Delta^1(\zeta)f)(s,y_1)g(s,y_1,y_2)\tq(t,x_1,x_2;s,y_1,y_2)  \d y_1\d y_2\\
&& = D_{x_{1}}^{N^1}D_{x_{2}}^{N^2} \int_{\R^{2d}}   (\Delta^1(\zeta)f)(s,y_1)(\Delta^2(\zeta)g)(s,y_1,y_2)\tq(t,x_1,x_2;s,y_1,y_2)  \d y_1\d y_2  \nonumber\\
&&\qquad + D_{x_{1}}^{N^1}D_{x_{2}}^{N^2}  \int_{\R^{2d}}  (\Delta^1(\zeta)f)(s,y_1) g(s,y_1,\zeta_2)  \tq(t,x_1,x_2;s,y_1,y_2)  \d y_1\d y_2. \nonumber
\end{eqnarray*}
By using differentiation under the integral sign Theorem, \eqref{argcentre} (since $N^2$ is positive) together with an integration by parts, the last term in the right hand side is equal to 0. Finally, assertion \eqref{c1centeresti2} follows from the same arguments. This concludes the proof of the Claim.
\end{proof}

\subsection{Representation and differentiation of the solution of the PDE \eqref{PDE2}}\label{rereso}
\begin{lemme}\label{repderivsol}
Suppose that assumptions \textbf{(HR)} hold, then, for all $x$ in $\R^{2d}$ and $t$ in $[0,T]$ the solution $u=(u_1,u_2)^*$ of the PDE \eqref{PDE2} can be written as
\begin{eqnarray}
u_i(t,x) &=& E\left[\int_t^T\phi_i(s,\tilde{X}^{t,x}_s) - (\mL-\tilde{\mL}^{t,\xi})u_i(s,\tilde{X}^{t,x}_s) \d s\right]\label{spit2}\\
&=&\int_t^T \int_{\R^{2d}} \phi_i(s,y)  \tilde{q}(t,x;s,y)\d y\d s\notag\\
&&- \int_t^T \int_{\R^{2d}} \frac{1}{2}{\rm Tr}\left[\Delta(\theta_{t,s}(\xi))a(s,y) D^2_{y_1}u_i(s,y)\right] \tilde{q}(t,x;s,y)\d y \d s \nonumber\\
&&- \int_t^T \int_{\R^{2d}} \Big\{\left[\Delta(\theta_{t,s}(\xi))F_1(s,y)\right]\cdot D_{y_1}u_i(s,y)\Big\} \tilde{q}(t,x;s,y)\d y\d s \nonumber\\
&&- \int_t^T \int_{\R^{2d}} \Big\{\left[\Delta(\theta_{t,s}(\xi))F_2(s,y)  - D_{1}F_2(s,\theta_{t,s}(\xi))\Delta^1(\theta_{t,s}(\xi))y_1\right]  \cdot D_{y_2}u_i(s,y) \Big\}\tilde{q}(t,x;s,y)\d y \d s\nonumber\\
&= :& \sum_{j=1}^4  \int_{t}^T \int_{\R^{2d}} \mH_i^j(s,y,\theta_{t,s}(\xi)) \tilde{q}(t,x;s,y) \d y\d s, \label{grgrgr1}\\
&=: & \sum_{j=1}^4 \int_t^T \mI^j_i(s,x)  \d s \label{repsol} 
\end{eqnarray}
for all $\xi$ in $\R^{2d}$. It is infinitely differentiable on $[0,T]\times \R^{2d}$.
\end{lemme}

\begin{remarque}
We emphasize that the solution ``$u$'' does not depend on the choice of the freezing data ``($\tau,\xi$)''. Above we chose to set $\tau=t$ when the solution $u$ is evaluated at time $t$. We keep this choice from now on.
\end{remarque}

\begin{proof}
Thanks to Lemma \ref{eeu} the PDE \eqref{PDE2} is well posed and can be rewritten as
\begin{equation}\label{casinterPDE}
\left\lbrace\begin{array}{ll}
\p_t u_i(t,x) + \tilde{\mL}^{\tau,\xi}u_i(t,x) =-(\mL -\tilde{\mL}^{\tau,\xi})u_i(t,x) +\phi_i(t,x),\ (t,x)\in [0,T)\times \R^{2d}\\ 
u_i(T,x)=0, \ i=1,2,
\end{array}\right.
\end{equation}
for all $(\tau,\xi)$ in $[0,T]\times \R^{2d}$, so that \eqref{spit2} follows from Feynman-Kac representation. Equality \eqref{repsol} follows from the definition of $\tilde{q}$ in Proposition \ref{estfundsol}. Next, given a positive $\epsilon $, we have for all $(t,x)$ in $[0,T]\times \R^{2d}$,

\begin{eqnarray}
u_i(t,x)  &=&\int_{t+\epsilon}^T\E\left[\phi_i(s,\tilde{X}^{t,x}_s) - (\mL-\tilde{\mL}^{\tau,\xi})u_i(s,\tilde{X}^{t,x}_s)\right]\d s \notag\\
&&+ \int_{t}^{t+\epsilon} \E\left[\phi_i(s,\tilde{X}^{t,x}_s) - (\mL-\tilde{\mL}^{\tau,\xi})u_i(s,\tilde{X}^{t,x}_s)\right]\d s\notag.
\end{eqnarray}
Under \textbf{(HR)}, the coefficients of $\mL$, $\tilde{\mL}^{\tau,\xi}$ and the functions $\phi_i$, $u_i$ are smooth (Lemma \ref{eeu}). From classical regularity results on the solution of the SDE \eqref{kolmogf} (see e.g. \cite{kunita_stochastic_1982}) we can deduce that the solution is infinitely differentiable under \textbf{(HR)}.
\end{proof}

We now derive a representation formula for the derivatives of $\mI_i^j,\ i=1,2,\ j=1,\ldots,4,$ that involve a differentiation w.r.t. the degenerate variable. This allows us to handle the singularity of the derivative of the kernel $\tilde{q}$  (see assertion (iv) in Proposition \ref{estfundsol}) as done in Subsection \ref{SZ4} (\emph{centering argument}). 

\begin{lemme}\label{center}\label{center2}
Let $t<T$ in $\R^+$. For all $i$ in $\{1,2\}$, all $(s,x)$ in $(t,T] \times \R^{2d}$ and all integer $n$ 

the terms $D_{x_1}^nD_{x_2}\mI^j_i(s,x)$, $j=1,\ldots,4,$ can be written as:
\begin{eqnarray}\label{centerH1}
D_{x_1}^nD_{x_2}\mI^1_i(s,x) &=&  - \int_{\R^{2d}}  \Delta^2(\theta_{t,s}(\xi))\phi_i(s,y) D^n_{x_1}D_{x_2}\tilde{q}(t,x;s,y) \d y 
\end{eqnarray}
and
\begin{eqnarray}\label{centerH2}
&&D_{x_1}^nD_{x_2}\mI_i^2(s,x)\\ 
&=&-\frac{1}{2}  \int_{\R^{2d}}\Bigg\{  {\rm Tr}\left[[\Delta^2(\theta_{t,s}(\xi))a(s,y)]D^2_{y_1}u_i(s,y)\right] \notag\\
&&\qquad - \sum_{l=1}^d\left[\frac{\p }{\p y_{1l}}a_{l.}(s,y)\right]. \left[\Delta^2(\theta_{t,s}(\xi))D_{y_1}u_i(s,y)\right]\Bigg\} D_{x_1}^nD_{x_2}\tilde{q}(t,x;s,y)\d y \nonumber\\
&& +\frac{1}{2}\sum_{l=1}^d  \int_{\R^{2d}} \bigg\{\left[\Delta^1(\theta_{t,s}(\xi))a_{l.}(s,y)\right]
\cdot \left[\Delta^2(\theta_{t,s}(\xi))D_{y_1}u_i(s,y)\right]\bigg\}  D_{x_1}^nD_{x_2} \left(\frac{\p}{\p {y_{1l}}}\tilde{q}(t,x;s,y)\right)  \d y ,\nonumber
\end{eqnarray}
where ``$a_{l.}$'' denotes the $l^{{\rm th}}$ line of the matrix $a$, and
\begin{eqnarray}\label{centerH3}
D_{x_1}^nD_{x_2}\mI_i^3(s,x) &=&- \int_{\R^{2d}}\Bigg\{ \left[\Delta^2(\theta_{t,s}(\xi))F_1(s,y)\right] \cdot D_{y_1}u_i(s,y)\\
&&\quad +\left[ \Delta^1(\theta_{t,s}(\xi))F_1(s,y)\right] \cdot \left[\Delta^2(\theta_{t,s}(\xi))D_{y_1}u_i(s,y)\right]\Bigg\}  D_{x_1}^n D_{x_2}\tilde{q}(t,x;s,y)\d y,\nonumber
\end{eqnarray}
and finally:
\begin{eqnarray}
D_{x_1}^nD_{x_2}\mI_i^4(s,x) &=& - \int_{\R^{2d}}\bigg\{\left[ \Delta^1(\theta_{t,s}(\xi))F_2(s,y) - D_{1}F_2(s,\theta_{t,s}(\xi))\Delta^1(\theta_{t,s}(\xi))y_1  \right] \label{centerH4}\\
&& \qquad \cdot \left[\Delta^2(\theta_{t,s}(\xi))D_{y_2}u_i(s,y) \right]+\left[\Delta^2(\theta_{t,s}(\xi))F_2(s,y)\right]\cdot   D_{y_2}u_i(s,y)\bigg\} \notag\\
&&\qquad \hphantom{space}\times D_{x_1}^nD_{x_2}\tilde{q}(t,x;s,y) \d y.\notag
\end{eqnarray}
\end{lemme}

\begin{proof}
Representation \eqref{centerH1} is a direct consequence of assertion \eqref{c1centeresti2} in  Claim \ref{c1center}. Next, we deal with \eqref{centerH2}. By using first the decomposition $\Delta=\Delta^1+\Delta^2$ and then by integrating by parts, we have
\begin{eqnarray}
\mI_i^2(s,x) &=& -\frac{1}{2} \int_{\R^{2d}} \Bigg\{{\rm Tr}\left[[\Delta^2(\theta_{t,s}(\xi))a(s,y)]D^2_{y_1}u_i(s,y)\right]\notag\\
&& \qquad-\sum_{l=1}^d\left[\frac{\p }{\p y_{1l}}\Delta^1(\theta_{t,s}(\xi))a_{l.}(s,y)\right] \cdot D_{y_1}u_i(s,y)\Bigg\} \tilde{q}(t,x;s,y) \d y \notag\\
&&  +\frac{1}{2} \sum_{l=1}^d  \int_{\R^{2d}} \Big\{\left[\Delta^1(\theta_{t,s}(\xi))a_{l.}(s,y)\right] \cdot D_{y_1}u_i(s,y)\Big\}\frac{\p}{\p y_{1l}}\tilde{q}(t,x;s,y) \d y .\notag
\end{eqnarray}
Note that, for all $l$ in $\{1,\cdots,d\}$, $[\p/\p y_{1l}]\Delta^1(\theta_{t,s}(\xi))a_{l.}(s,y_1,\theta_{t,s}^2(\xi))=[\p/\p y_{1l}]a_{l.}(s,y_1,\theta_{t,s}^2(\xi))$. We conclude by differentiating and then by applying assertion \eqref{c1centeresti3} in Claim \ref{c1center} with $f=a$ and $g=D_{1}u_i$. This gives \eqref{centerH2}.

By using again first the decomposition $\Delta=\Delta^1+\Delta^2$, assertions \eqref{centerH3} and \eqref{centerH4} are immediate consequences of assertion \eqref{c1centeresti3} in Claim \ref{c1center}. This concludes the proof of Lemma \ref{center}.
\end{proof}

We conclude this section with the following definition:

\begin{definition}\label{integrands}
\textbf{The centered integrands.} When differentiating at least once the terms $\mI^j_i$, $j=1,\ldots,4,$ w.r.t. the degenerate variable ``$x_2$'', Lemma \ref{center} allows us to identify the $\mH_i^j$, $j=1,3,4,$ defined in \eqref{grgrgr1} with the integrand of \eqref{centerH1}, \eqref{centerH3}, and \eqref{centerH4} respectively. In the same spirit we use in that case the notation $\mH_i^2$ to denote both the sum of the two integrands appearing in the representation \eqref{centerH2} of $D_{x_1}^nD_{x_2}\mI_i^2$ and the integrand defined in \eqref{grgrgr1}. 

When we identify the $\mH^j_i$, $j=1,\ldots,4,$ with the integrand of \eqref{centerH1}, \eqref{centerH2}, \eqref{centerH3} and \eqref{centerH4} respectively, we will call them the \emph{centered integrands}.
\end{definition}

%
%
%

\section{Proof of Proposition \ref{EE}}\label{PBFGH}

\subsection{From intermediate gradient estimates to Proposition \ref{EE}}\label{FPTULE}

We here give the intermediate estimates that allow to prove Proposition \ref{EE}. In the following, $u_i$ denotes the $i^{th}$ component of the solution $u=(u_1,u_2)^*$ of the linear system of PDE (\ref{PDE2}). The following arguments and Lemmas hold for $i=1,2$.

Since the representation \eqref{repsol} of each $u_i$ involves its derivatives, we prove Proposition \ref{EE} by using a \emph{circular argument} (see \eqref{circular}). We recall that afor ny function from $[0,T]\times \R^d \times \R^d$ we denote by $D_1$ (resp. $D_2$) the derivative with respect to the first (resp. second) $d$-dimensional space component. We first show that

\begin{lemme}\label{ouf}
Suppose that assumptions \textbf{(HR)} hold. There exist a positive $\mc{T}_{\ref{ouf}}$, two positive numbers $\delta_{\ref{ouf}}$ and $\bar{\delta}_{\ref{ouf}}$ and a positive constant $C$, depending only on known parameters in \textbf{(H)}, such that
\begin{enumerate}[(i)]
\item $\displaystyle \left\|D^2_{1}u_i\right\|_{\infty} \leq T^{\delta_{\ref{ouf}}} C \left( 1+ \left\|D_{2}u_i\right\|_{\infty}\right),$\label{urgu1}
\item $\displaystyle \left\|D_{1}u_i\right\|_{\infty}\leq T^{\bar{\delta}_{\ref{ouf}}} C\left(1+\left\|D_{2}u_i\right\|_{\infty}\right),
$\label{urgu2}
\end{enumerate}
for all $T$ less than $\mc{T}_{\ref{ouf}}$. We recall that $D_1$ denotes the derivative with respect to the first $d$-dimensional space component
\end{lemme}

Then, we estimate the gradients that involve the derivatives w.r.t. the degenerate variable ``$D_2$''. To this aim, we differentiate the representation \eqref{repsol} and we estimate it. Thanks to Proposition \ref{estfundsol}, we know that this differentiation generates the worst singularity in the time space integrals. As shown in Subsection \ref{SZ4}, we can use the regularity of the coefficients assumed in \textbf{(H)} in order to smooth this singularity. But in this more general case, we also have to use the regularity of the solution itself. Notably, we need to estimate the H\"{o}lder regularity of $D_{2}u_i$. Hence, we prove the following sort of H\"{o}lder estimate on $D_{2}u_i$:

\begin{lemme}\label{holdregderivsecdegen}
Suppose that assumptions \textbf{(HR)} hold and let 
\begin{equation}\label{defdeM}
M(D_{2}u_i,T) := \sup_{w_1,w_2 \neq w'_2 \in \R^{d},\ t\in [0,T]} \frac{|D_{2}u_i(t,w_1,w_2)-D_{2}u_i(t,w_1,w'_2)|}{|w_2-w'_2|^{\gamma/3}+|w_2-w'_2|^{\beta_2^2}+|w_2-w'_2|^{\beta_1^2}+|w_2-w'_2|},
\end{equation}
for some positive number $\gamma$. There exist a positive $\mc{T}_{\ref{holdregderivsecdegen}}$, a positive constant $C$ and a positive number $\delta_{\ref{holdregderivsecdegen}}$, depending only on known parameters in \textbf{(H)}, such that for all positive $\gamma$ strictly less than $ 3\inf\{\beta_1^2,\beta_2^2\}-1$, 
\begin{equation*}
M(D_{2}u_i,T)\leq C T^{\delta_{\ref{holdregderivsecdegen}}} \left(1+\left\|D_{2}u_i\right\|_{\infty} +\left\|D_{1}D_{2}u_i\right\|_{\infty} \right),
\end{equation*}
for all $T$ less than $\mc{T}_{\ref{holdregderivsecdegen}}$.
\end{lemme}

This allows us to obtain the following estimates

\begin{lemme}\label{supderivsec}
Suppose that assumptions \textbf{(HR)} hold and let $n$ in $\{0,1\}$. There exist a positive $\mc{T}_{\ref{supderivsec}}(n)$, a positive number $\delta_{\ref{supderivsec}}(n)$ and a positive constant $C(n)$, depending only on known parameters in \textbf{(H)} and $n$, such that:
\begin{equation*}
\left\|D_{1}^nD_{2}u_i\right\|_{\infty} \leq C(n) T^{\delta_{\ref{supderivsec}}(n)},
\end{equation*}
for all $T$ less than $\mc{T}_{\ref{supderivsec}}(n)$.
\end{lemme}

Finally, we can deduce Proposition \ref{EE} from Lemmas \ref{ouf} and \ref{supderivsec}. \qed

\subsection{Intermediate gradient estimates: proofs}
\subsubsection{Proof of Lemma \ref{ouf} }\label{SNDPDSPRDD}
We first show assertion \eqref{urgu1}.  Let $(t,x)$ belong to $[0,T]\times \R^{2d}$, from the representation \eqref{repsol} in Lemma \ref{repderivsol} we have:
\begin{equation}\label{iyteri}
D^2_{x_1}u_i(t,x) = \sum_{j=1}^{4} D^2_{x_1}\int_t^T \mI_i^j(s,x) \d s.
\end{equation} 
As done in the subsection \ref{SZ4} we have to obtain a suitable bound on the $D^2_{1}\mI_i^j$ in order to estimate $D^2_{1}u_i$. Thanks to assertion \eqref{c1centeresti1} in Claim \ref{c1center} applied on the integrand of $D^2_{x_1}\mI_i^1$ we have for all $s$ in $[t,T]$:
\begin{eqnarray}
\sum_{j=1}^{4} D^2_{x_1} \mI_i^j(s,x)  &=&  \int_{\R^{2d}} \Delta(\theta_{t,s}(\xi))\phi_i(s,y) \left[D^2_{x_1}\tilde{q}(t,x;s,y)\right]  \d y \label{notag}\\
&&  - \frac{1}{2} \int_{\R^{2d}} {\rm Tr}\left[ \Delta(\theta_{t,s}(\xi))a(s,y) D^2_{y_1}u_i(s,y)\right] \left[D^2_{x_1}\tilde{q}(t,x;s,y)\right] \d y \nonumber\\
&&  -  \int_{\R^{2d}} \left\{\left[\Delta(\theta_{t,s}(\xi))F_1(s,y)\right] \cdot  D_{y_1}u_i(s,y)\right\} \left[D^2_{x_1}\tilde{q}(t,x;s,y)\right] \d y\nonumber\\
&&  - \int_{\R^{2d}} \Big\{\left[ \Delta(\theta_{t,s}(\xi))F_2(s,y) - D_{1}F_2(s,\theta_{t,s}^1(\xi),\theta_{t,s}^2(\xi)) \Delta^1(\theta_{t,s}(\xi))y_1  \right]\notag\\
&& \quad \cdot \left[D_{y_2}u_i(s,y) \right] \Big\}\left[D^2_{x_1}\tilde{q}(t,x;s,y)\right]  \d y .\nonumber
\end{eqnarray}

Note that by using the fact that $\Delta=\Delta^1+\Delta^2$, a Taylor expansion of order 0 with integrable remainder of the mapping $\theta_{t,s}^1(\xi) \mapsto F_2(\cdot,\theta_{t,s}^1(\xi),\cdot)$ around $y_1$ in $\R^d$ and \textbf{(H3-a)}, we have that:
\begin{eqnarray}
\forall (s,y_2) \in (t,T]\times \R^d, && \left|\left[\Delta(\theta_{t,s}(\xi))F_2(s,y)  - D_{1}F_2(s,\theta_{t,s}(\xi)) \left(\Delta^1(\theta_{t,s}(\xi))y_1\right) \right] \cdot D_{y_2}u_i(s,y)\right|\notag\\
&&\quad \leq C\left\|D_{2}u_i\right\|_{\infty} \left(|\Delta^2(\theta_{t,s}(\xi))y_2|^{\beta_2^2} + |\Delta^1(\theta_{t,s}(\xi))y_1|^{1+\alpha^1}\right).\label{taylorexp}
\end{eqnarray}

From Proposition \ref{estfundsol}, we know that for all $s$ in $(t,T]$ and $y$ in $\R^{2d}$, $|D^2_{x_1}\tilde{q}(t,x;s,y)| \leq C' (s-t)^{-1} \hat{q}_{c}(t,x;s,y)$. By plugging this estimate in \eqref{notag}, together with the regularity of the coefficients assumed in \textbf{(H)} and \eqref{taylorexp}, we obtain that

\begin{eqnarray*}
\left|\sum_{j=1}^{4} D^2_{x_1} \mI_i^j(s,x) \right| &\leq & C''  (s-t)^{-1} \int_{\R^{2d}}  \Bigg\{\sum_{j=1}^2 \Bigg\{(s-t)^{(j-1/2)\beta_i^j}\left|\frac{\Delta^j(\theta_{t,s}(\xi))y_j}{(s-t)^{(j-1/2)}}\right|^{\beta_i^j}+ \left\|D^2_{1}u_i\right\|_{\infty}  \\
&&\quad   \times (s-t)^{(j-1/2)}  \left|\frac{\Delta^j(\theta_{t,s}(\xi))y_j}{(s-t)^{(j-1/2)}}\right| +  \left\|D_{1}u_i\right\|_{\infty}  (s-t)^{(j-1/2)\beta^j_1}\left|\frac{\Delta^j(\theta_{t,s}(\xi))y_j}{(s-t)^{(j-1/2)}}\right|^{\beta^j_1} \Bigg\}\\
&&  +\left\|D_{2}u_i\right\|_{\infty}  \Bigg[(s-t)^{3\beta_2^2/2}  \left|\frac{\Delta^2(\theta_{t,s}(\xi))y_2}{(s-t)^{3/2}}\right|^{\beta_2^2} +  (s-t)^{(1+\alpha^1)/2}\\
&&\quad \times   \left|\frac{\Delta^1(\theta_{t,s}(\xi))y_1}{(s-t)^{1/2}}\right|^{1+\alpha^1} \Bigg] \Bigg\}  \hat{q}_{c}(t,x;s,y)\d y.
\end{eqnarray*}
Set $\xi=x$, we now use the Gaussian off-diagonal decay in $\hat{q}_c$ (see the computations in Subsection \ref{SZ4} for more details):
\begin{eqnarray}\label{gaussdecay}
\forall \kappa>0, \exists  \underbar{C}, \underbar{c}>0 \text{ s.t. } \forall (s,y) \in (t,T] \times \R^{2d}:\ |\Delta^l(\theta_{t,s}(x))y_l|^\kappa \hat{q}_c(t,x;s,y) \leq \underbar{C} (s-t)^{(l-1/2)\kappa} \hat{q}_{\underbar{c}}(t,x;s,y)
\end{eqnarray}
for $l=1,2$ and where $\underbar{C}$ and $\underbar{c}$ depend only on known parameters in \textbf{(H)} and on $\kappa$. Hence, by integrating w.r.t the space variable we have:

\begin{eqnarray*}
&&\left|\sum_{j=1}^{4} D^2_{x_1} \mI_i^j(s,x)\right|\\
&&\leq C''' \Bigg\{\sum_{j=1}^2 \Bigg\{ (s-t)^{-1+(j-1/2)\beta_i^j}+  \left((s-t)^{-1/2}+(s-t)^{1/2} \right)\left\|D^2_{1}u_i\right\|_{\infty} \\
&&\quad   +  (s-t)^{(j-1/2)\beta^j_2-1}  \left\|D_{1}u_i\right\|_{\infty} \Bigg\} +   \Big((s-t)^{-1+3\beta_1^2/2} +(s-t)^{(\alpha^1-1)/2}\Big) \left\|D_{2}u_i\right\|_{\infty} \Bigg\} .
\end{eqnarray*}
So that we can inverse the differentiation and integral operators in \eqref{iyteri} and then deduce that

\begin{eqnarray*}
\left\|D^2_{1}u_i\right\|_{\infty} &\leq& C''' \bigg\{ \left( T^{3\beta_2^2/2} + T^{(1+\alpha^1)/2} \right)\left\|D_{2}u_i\right\|_{\infty}+ \left( T^{\beta_1^1/2} + T^{3\beta_1^2/2} \right) \left\|D_{1}u_i\right\|_{\infty} \\
&&+  \left(T^{1/2}+T^{3/2}\right) \left\|D^2_{1}u_i\right\|_{\infty} + T^{\beta_i^1/2} + T^{3\beta_i^2/2}\bigg\}.
\end{eqnarray*}
Then, by setting $\mc{T}^{\eqref{urgu1}}_{\ref{ouf}} = \sup\left\{T>0,\ {\rm such\ that}\ C'''(T^{1/2}+T^{3/2}) \leq 1/2\right\}$  we deduce the assertion \eqref{urgu1} for all $T$ less than $\mc{T}^{\eqref{urgu1}}_{\ref{ouf}}$ from a \emph{circular argument} (see \eqref{circular}). The proof of the second statement \eqref{urgu2} can be done by the same arguments. This concludes the proof of Lemma \ref{ouf}.\qed

\subsubsection{Proof of Lemma \ref{supderivsec}}\label{ESNFG}

 Let $(t,x)$ belongs to $[0,T]\times \R^{2d}$. We have from representation \eqref{repsol} in Lemma \ref{repderivsol} that, for all $n$ in $\{0,1\}$ and all positive $\epsilon$
\begin{equation}\label{lareprororo}
D_{x_1}^{n}D_{x_2}u_i(t,x) = \sum_{j=1}^{4} D_{x_1}^{n} D_{x_2} \int_t^T \mI_i^j(s,x) \d s.
\end{equation}
Again, we look for a suitable bound on the $D_{x_1}^{n} D_{x_2} \mI_i^j$ in order to estimate $D_{x_1}^{n} D_{x_2}u_i$.
Thanks to representation \eqref{grgrgr1}, we know that the derivatives $D_{x_1}^{n} D_{x_2}\mI_i^j$ can be written as the integral of some function against the derivative of the degenerate Gaussian kernel $\tilde{q}$. Since from Proposition \ref{estfundsol} we have that for all $s$ in $(t,T]$ and $y$ in $\R^{2d}$
\begin{equation}\label{estigeussinterrrr}
|D_{x_1}^nD_{x_2}\tilde{q}(t,x;s,y)| \leq C (s-t)^{-(3+n)/2} \hat{q}_{c}(t,x;s,y),
\end{equation} 
we use the regularity of the coefficients together with the regularity of the solution $u_i$ and its derivatives to smooth the time singularity appearing in \eqref{estigeussinterrrr}, as done in the previous subsection and in Subsection \ref{SZ4}.\\

Hence, we first give a bound on each \emph{centered integrand} $\mH_i^j$ of $D_{x_1}^nD_{x_2}\mI_i^j,$ $j=1,\cdots,4,$ given by Definition \ref{integrands}. For all $s$ in $(t,T]$, $y$ in $\R^{2d}$ we have, from the regularity of the coefficients assumed in \textbf{(H)}, that:
\begin{equation}\label{estiF1}
|\mH_i^1(s,y,\theta_{t,s}(\xi))| \leq C \left|\Delta^2(\theta_{t,s}(\xi))y_2\right|^{\beta_i^2}.
\end{equation}
Then, we recall that from Mean Value Theorem (MVT) we have
\begin{equation}\label{labelle1}
\left|\Delta^2(\theta_{t,s}(\xi))D_{y_1}u_i(s,y)\right| \leq \left\|D_{1}D_{2}u_i\right\|_{\infty}|\Delta^2(\theta_{t,s}(\xi))y_2|,
\end{equation}
so that 
\begin{eqnarray}\label{estiF2}
|\mH_i^2(s,y,\theta_{t,s}(\xi))|  & \leq & C'\bigg\{ \left\|D^2_{1}u_i\right\|_{\infty} \left|\Delta^2(\theta_{t,s}(\xi))y_2\right|^{} + \left\|D_{1}D_{2}u_i\right\|_{\infty}\notag\\
&& \times  \bigg( \left|\Delta^2(\theta_{t,s}(\xi))y_2\right|^{}+ \left|\Delta^1(\theta_{t,s}(\xi))y_1\right|^{}\left|\Delta^2(\theta_{t,s}(\xi))y_2\right|^{}\bigg)\bigg\}.
\end{eqnarray}
By the same way, we get that
\begin{eqnarray}\label{estiF3}
|\mH_i^3(s,y,\theta_{t,s}(\xi))| & \leq & C''\Big\{ \left\|D_{1}u_i\right\|_{\infty} \left|\Delta^2(\theta_{t,s}(\xi))y_2\right|^{\beta_1^2} \\
&&+ \left\|D_{1}D_{2}u_i\right\|_{\infty} \left|\Delta^1(\theta_{t,s}(\xi))y_1\right|^{\beta_1^1}\left|\Delta^2(\theta_{t,s}(\xi))y_2\right|^{}\Big\}.\notag
\end{eqnarray}
And finally, since we have from Lemma \ref{holdregderivsecdegen}
\begin{eqnarray}\label{labelle2}
&&\left|\Delta^2(\theta_{t,s}(\xi))D_{y_2}u_i(s,y)\right| \\
&& \leq  M(D_{2}u_i,T) \Bigg(|\Delta^2(\theta_{t,s}(\xi))y_2|^{\gamma/3} + |\Delta^2(\theta_{t,s}(\xi))y_2|^{\beta_1^2}  +|\Delta^2(\theta_{t,s}(\xi))y_2|^{\beta_2^2} + |\Delta^2(\theta_{t,s}(\xi))y_2|^{}\Bigg),\notag
\end{eqnarray}
where $M(D_{2}u_i,T) $ is defined by \eqref{defdeM} we have from the regularity of the coefficients assumed in \textbf{(H)} and \eqref{taylorexp} that
\begin{eqnarray}\label{estiF4}
|\mH_i^4(s,y,\theta_{t,s}(\xi))|  &\leq & C''' \left\|D_{2}u_i\right\|_{\infty} \left|\Delta^2(\theta_{t,s}(\xi))y_2\right|^{\beta_2^2} + M(D_{2}u_i,T) |\Delta^1(\theta_{t,s}(\xi))y_1|^{1+ \alpha^1} \\
&&\times \Bigg(|\Delta^2(\theta_{t,s}(\xi))y_2|^{\gamma/3} + |\Delta^2(\theta_{t,s}(\xi))y_2|^{\beta_1^2}  +|\Delta^2(\theta_{t,s}(\xi))y_2|^{\beta_2^2} + |\Delta^2(\theta_{t,s}(\xi))y_2|^{}\Bigg).\notag
\end{eqnarray}

Finally, let us recall from Proposition \ref{estfundsol} that  for all $(s,y)$ in $(t,T] \times \R^{2d}$,
\begin{equation}\label{estigeussinterrrr2r}
\left|D_{x_1}^nD_{x_2}\frac{\p}{\p y_{1l}}\tilde{q}(t,x;s,y)\right| \leq C'''' (s-t)^{-(4+n)/2} \hat{q}_{c'}(t,x;s,y),
\end{equation} 
for all $l$ in $\{1,\ldots,d\}$.

Now, we can plug together some of the above estimates in the corresponding $D_{x_1}^{n}D_{x_2}\mI_i^j$, $j=1,\ldots,4,$ defined in Lemma \ref{center}. By using  \eqref{estigeussinterrrr} with \eqref{estiF1} (resp. \eqref{estiF3} and \eqref{estiF4}) in \eqref{centerH1} (resp. \eqref{centerH3} and \eqref{centerH4}), estimates \eqref{estigeussinterrrr}  and \eqref{estigeussinterrrr2r} with \eqref{estiF2} in \eqref{centerH2}, by letting next $\xi=x$ and by using \eqref{gaussdecay} in all these terms, we can deduce that for all $s$ in $(t,T]$
%
\begin{eqnarray}
&&\left|\sum_{j=1}^{4} D_{x_1}^{n} D_{x_2} \mI_i^j(s,x)  \right| \label{zereff}\\
&& \leq c'\int_{\R^{2d}}\Bigg\{ (s-t)^{(3(\beta_i^2-1)-n)/2} + (s-t)^{-n/2}\left( \left\|D^2_{1}u_i\right\|_{\infty} + \left\|D_{1}D_{2}u_i\right\|_{\infty}\right)\notag \\
&&  + \|D_{1}u_i\|_{\infty} (s-t)^{(3(\beta_1^2-1)-n)/2}   +  \left\|D_{1}D_{2} u_i\right\|_{\infty}  (s-t)^{(\beta_1^1-n)/2} + M(D_{2}u_i,T)   (s-t)^{-1+(\alpha^1-n)/2}  \notag\\
&& \quad \times \left((s-t)^{\gamma/2}+(s-t)^{\beta_1^2/2}+(s-t)^{\beta_2^2/2}+(s-t)^{}\right) + \left\|D_{2}u_i\right\|_{\infty}  (s-t)^{(3(\beta_2^2-1)-n)/2}\Bigg\} \notag\\
&&\qquad \times \hat{q}_{c''}(t,x;s,y)\d y.\notag
\end{eqnarray}

Then, we have from Lemma \ref{holdregderivsecdegen} that for all $T$ less than $\mc{T}_{\ref{holdregderivsecdegen}}$:
\begin{equation}\label{interMMMMMM}
\forall\gamma\ {\rm s.t.}\  0<\gamma< 3\inf\{\beta_1^2,\beta_2^2\}-1,\quad  M(D_{2}u_i,T)\leq c''' \mc{T}_{\ref{holdregderivsecdegen}} \left(1+\left\|D_{2}u_i\right\|_{\infty} +\left\|D_{1}D_{2}u_i\right\|_{\infty} \right).
\end{equation}
So that we can inverse the differentiation and integral operators in \eqref{lareprororo}. Hence, by plugging \eqref{interMMMMMM} in \eqref{zereff} and by choosing $n$ equal to 0, after integrating in time,  we deduce from a \emph{circular argument} (as described in \eqref{circular}) that there exist two positive numbers $\mathcal{T}'_{\ref{supderivsec}}(0)$, $\delta_{\ref{supderivsec}}(0)$ and a positive constant $\bar{C}$, depending only on known parameters in (\textbf{H}), such that
\begin{eqnarray}\label{interterrible}
\left\|D_{2}u_i\right\|_{\infty} &\leq &\bar{C}T^{\delta_{\ref{supderivsec}}(0)}\left(1+\left\|D_{1}D_{2}u_i\right\|_{\infty}\right),
\end{eqnarray}
for all $T$ less than $\mathcal{T}'_{\ref{supderivsec}}(0)$. Next, by letting $n$ be equal to $1$ in \eqref{zereff} (recall that Lemma \ref{holdregderivsecdegen} hold for all $0<\gamma< 3\inf\{\beta_1^2,\beta_2^2\}-1$ and that $3\inf\{\beta_1^2,\beta_2^2\}>1$ so that all the time-singularity are integrable), by using the same arguments as in the case $n=0$ together with \eqref{interterrible}, we can show that there exist two positive numbers $\mathcal{T}_{\ref{supderivsec}}(1)$, $\delta_{\ref{supderivsec}}(1)$ and a positive constant $\bar{C}'$, depending only on known parameters in (\textbf{H}), such that
\begin{eqnarray*}
\left\|D_{1}D_{2}u_i \right\|_{\infty} &\leq &\bar{C}'T^{\delta_{\ref{supderivsec}}(1)},
\end{eqnarray*}
for all $T$ less than $\mc{T}_{\ref{supderivsec}}(1)$. This conclude the proof for $n=1$. The case $n=0$ follows from plugging this estimate in \eqref{interterrible}.\\

\begin{remarque}\label{remarqueee} Note that the estimate on supremum norm of $D_2u_i$ (\emph{i.e.} when $n=0$) could be obtain wihtout Lemma \ref{holdregderivsecdegen}, by bounding the left hand side of \eqref{labelle2} by the supremum norm of $D_2u_i$.
\end{remarque}

\subsubsection{Proof of Lemma \ref{holdregderivsecdegen}}\label{HRDD}
From \eqref{repsol} in Lemma \ref{repderivsol}, Lemma \ref{supderivsec} and Remark \ref{remarqueee}, for all $(t,x_1)$ in $[0,T]\times \R^{d}$ and $x_2 \neq z_2$ in $\R^{d}$ we have:
\begin{eqnarray}\label{sumterms}
&&\left|D_{2}u_i(t,x_1,x_2) -  D_{2}u_i(t,x_1,z_2)\right| \\
&&\leq \left| \sum_{j=1}^4  \int_t^T D_{x_2}\mI^j_i(s,x_1,x_2) - D_{z_2}\mI^j_i(s,x_1,z_2) \d s\right|.\notag
\end{eqnarray}
We recall that the $\mI_i^j,\ j=1,\cdots,4,$ depend on the freezing point $\xi=(\xi_1,\xi_2)^*$ of the process which started from $(x_1,x_2)$ and $(x_1,z_2)$  at time $t$. Here, we choose the same freezing point ``$\xi$'' for the two processes (with different initial conditions). As done before, we  have to bound the dif

Then, we split the time interval w.r.t. the characteristic scale of the second component of the system \eqref{LS} in order to study the perturbation on each interval. Hence, we set $\mathcal{S}= \{s\in (t,T]\ {\rm s.t.}\ |x_2-z_2|<(s-t)^{3/2}\}$ and $\mathcal{S}^c= \{s\in (t,T]\ {\rm s.t.}\ |x_2-z_2|\geq(s-t)^{3/2}\}$. We have for all $s$ in $(t,T]$:

\begin{eqnarray}
&&\left| \sum_{j=1}^4  D_{x_2}\mI^j_i(s,x_1,x_2)- D_{z_2}\mI^j_i(s,x_1,z_2) \right|\notag\\
&& =\sum_{j=1}^4\ind \int_{\R^{2d}} \bigg\{\mH_i^j(s,y_1,y_2,\theta_{t,s}^1(\xi),\theta_{t,s}^2(\xi)) \notag\\
&& \qquad \bigg(D_{x_2}\tilde{q}(t,x_1,x_2;s,y_1,y_2) - D_{z_2}\tilde{q}(t,x_1,z_2;s,y_1,y_2)\bigg) \bigg\}\d y_1\d y_2\notag\\
&& \quad + \sum_{j=1}^4\indc \int_{\R^{2d}} \bigg\{\mH_i^j(s,y_1,y_2,\theta_{t,s}^1(\xi),\theta_{t,s}^2(\xi))\notag\\
&&\qquad \bigg(D_{x_2}\tilde{q}(t,x_1,x_2;s,y_1,y_2) - D_{z_2}\tilde{q}(t,x_1,z_2;s,y_1,y_2)\bigg) \bigg\}\d y_1\d y_2 \notag\\
&& :=  \mP_s(\mathcal{S}) +  \mP_s(\mathcal{S}^c),\label{omega}
\end{eqnarray}
where the $\mH_i^j,\ j=1,\cdots,4,$ are given by \eqref{grgrgr1}. We now bound the first and the second sum in the right hand side of the last equality above.\\

\textbf{Estimation of $\mP_s(\mathcal{S}) $ on $(t,T]$.} As a first step, we bound the first sum in the right hand side of \eqref{omega}. At the end of this part, it is proven that:
\begin{claim}\label{estifundasolholder}
For all $s$ in $\mathcal{S}$, $y_1,y_2$ in $\R^d$, the following inequality holds:
\begin{eqnarray}
&&\left|D_{x_2}\tilde{q}(t,x_1,x_2;s,y_1,y_2)-D_{z_2}\tilde{q}(t,x_1,z_2;s,y_1,y_2)\right|\nonumber\\
&& \quad \leq C(s-t)^{-(3+\gamma)/2}\hat{q}_{c}(t,x_1,x_2;s,y_1,y_2)|x_2-z_2|^{\gamma/3},\label{vmtdt}
\end{eqnarray}
where $c$ and $C$ depend only on known parameters in \textbf{(H)} and where $0<\gamma<3$. Moreover
\begin{eqnarray}
&&\left|\frac{\p}{\p y_{1 l}}D_{x_2}\tilde{q}(t,x_1,x_2;s,y_1,y_2)-\frac{\p}{\p y_{1 l}}D_{z_2}\tilde{q}(t,x_1,z_2;s,y_1,y_2)\right|\nonumber\\
&& \quad \leq C'(s-t)^{-(3+\gamma+1)/2}\hat{q}_{c'}(t,x_1,x_2;s,y_1,y_2)|x_2-z_2|^{\gamma/3},\label{vmtdtbis}
\end{eqnarray}
for any $l=1,\ldots,d$, where $c'$ and $C'$ depend only on known parameters in \textbf{(H)} and where $0<\gamma<3$.
\end{claim}


On the other hand, Lemma \ref{center} allows to choose the $\mH_i^j$, $j=1\ldots,4,$ appearing in $\mP_\epsilon(\mathcal{S})$ as the \emph{centered integrands} (see Definition  \ref{integrands}). So that this term is similar to the term studied in the proof of Lemma \ref{supderivsec}. In fact, the only difference with the terms in the proof of Lemma \ref{supderivsec} (when $n=0$) is that the integration is now against the perturbed differentiated kernel defined by the left hand side of \eqref{vmtdt} or \eqref{vmtdtbis}\footnote{And that the time integration is done on the set $\mathcal{S}$.}: in comparison with \eqref{estigeussinterrrr} and \eqref{estigeussinterrrr2r},  we ``loose'' a $\gamma/2$ in the time-singularity in the estimation. 

Hence, we can use the same arguments as the ones of Lemma \ref{supderivsec} to bound $\mP_s(\mathcal{S})$: by using estimates \eqref{estiF1}, \eqref{estiF2}, \eqref{estiF3}, \eqref{estiF4} together with \eqref{vmtdt} or \eqref{vmtdtbis} in $s(\mathcal{S})$ and next by letting $\xi=x$ together with \eqref{gaussdecay} we deduce that

\begin{eqnarray*}
\left|\mP_s(\mathcal{S})\right|  & \leq & C'   (s-t)^{-\gamma/2}\bigg\{ (s-t)^{-3(1-\beta_i^2)/2} + \left\|D^2_{1}u_i\right\|_{\infty} +\left\|D_{1}D_{2}u_i\right\|_{\infty} \notag \\
&&  + \| D_{1}u_i\|_{\infty} (s-t)^{-3(1-\beta_1^2)/2} +\|D_{1}D_{2}u_i\|_{\infty} (s-t)^{\beta_1^1/2} \notag\\ 
&& + \left\|D_{2}u_i\right\|_{\infty} (s-t)^{-3(1-\beta^2_2)/2} +M(D_{2}u_i,T) (s-t)^{-1+\alpha^1/2}   \notag\\
&& \quad \times\Big((s-t)^{\gamma/2}+ (s-t)^{3\beta_1^2/2} +  (s-t)^{3\beta_2^2/2}  +(s-t)^{3/2}\Big)\bigg\}  |x_2-z_2|^{\gamma/3},\notag
\end{eqnarray*}
for all $\gamma <3\inf(\beta_1^2,\beta_2^2)-1$. By integrating in space and by using Lemma \ref{ouf}, we obtain that there exists a positive number $\delta_{\ref{holdregderivsecdegen}}^1$, depending only on known parameters in \textbf{(H)},  such that:

\begin{eqnarray}\label{la1}
\int_t^T \left| \mP_s(\mathcal{S}) \right|\d s
&\leq &  C''T^{\delta_{\ref{holdregderivsecdegen}}^1}\left( M(D_{2}u_i,T)+ ||D_{2}u_i||_{\infty} + ||D_{1}D_{2}u_i||_{\infty} + 1\right)\label{la2}\\
&&\quad \times\left(|x_2-z_2|^{\gamma/3}+|x_2-z_2|^{\beta_2^2}+|x_2-z_2|^{\beta_1^2}+|x_2-z_2|\right), \notag
\end{eqnarray}
for all $\gamma < 3\inf (\beta_1^2,\beta_2^2)-1$.\\

\textbf{Estimation of $\mP_\epsilon(\mathcal{S}^c) $.} As a second step, we bound the sum $\mP_s(\mathcal{S}^c)$ on $(t,T]$ in \eqref{omega}. Note that
\begin{eqnarray}
\mP_s(\mathcal{S}^c) &=&\sum_{j=1}^4\indc \int_{\R^{2d}} \mH_i^j(s,y_1,y_2,\theta_{t,s}^1(\xi),\theta_{t,s}^2(\xi)) D_{x_2}\tilde{q}(t,x_1,x_2;s,y_1,y_2) \d y_1\d y_2\label{omegaga1}\\
& -& \sum_{j=1}^4\indc \int_{\R^{2d}} \mH_i^j(s,y_1,y_2,\theta_{t,s}^1(\xi),\theta_{t,s}^2(\xi)) D_{z_2}\tilde{q}(t,x_1,z_2;s,y_1,y_2)\d y_1\d y_2\label{omegaga2}
\end{eqnarray}
and that for all $s$ in $\mathcal{S}^c$ we have: 
\begin{equation}\label{inko}
1 \leq (s-t)^{-\gamma/2} \left|x_2-z_2\right|^{\gamma/3}.
\end{equation}

On the one hand we bound the right hand side of \eqref{omegaga1}. Again, thanks to Lemma \ref{center}, we can identify the $\mH_i^j$, $j=1\ldots,4,$ appearing in this term as the \emph{centered integrands} (see Definition  \ref{integrands}). And again, by proceeding exactly as in the proof of Lemma \ref{supderivsec} when $n=0$ (the restriction of the time integration on the set $\mathcal{S}^c$ is not a problem to do so), we obtain that this term is bounded by the right hand side of \eqref{zereff}. By using  \eqref{inko}, we deduce that

\begin{eqnarray}\label{hudgeinter}
&&\left|\sum_{j=1}^4\indc \int_{\R^{2d}} \mH_i^j(s,y_1,y_2,\theta_{t,s}^1(\xi),\theta_{t,s}^2(\xi)) D_{x_2}\tilde{q}(t,x_1,x_2;s,y_1,y_2) \d y_1\d y_2\right|  \\
&&\leq  C (s-t)^{-\gamma/2}  \Bigg\{ (s-t)^{-3(1-\beta_1^2)/2} + \left\|D^2_{1}u_i\right\|_{\infty} +\left\|D_{1}D_{2}u_i\right\|_{\infty} + \| D_{1}u_i\|_{\infty} (s-t)^{-3(1-\beta_1^2)/2}  \notag \\
&& \quad  +\|D_{1}D_{2}u_i\|_{\infty} (s-t)^{-\beta_1^1/2}  + \left\|D_{2}u_i\right\|_{\infty} (s-t)^{3(\beta^2_2-1)/2} + M(D_{2}u_i,T) (s-t)^{-1+\alpha^1/2}\notag\\
&&\qquad \times  \bigg( (s-t)^{\gamma/2}+(s-t)^{3\beta_1^2/2}+(s-t)^{3\beta_2^2/2}+(s-t)^{3/2} \bigg) \Bigg\} |x_2-z_2|^{\gamma/3}.\notag
\end{eqnarray}

%

On the other hand, we have to deal with the  term \eqref{omegaga2}. The crucial point here is that the frozen transition density is evaluated at point $z_2$ but frozen along the transport of $x_2$ (since we already chose $\xi=x$). This is because we took the same freezing point for the two solutions with different initial conditions. Hence, we have to re-center carefully each \emph{integrand} $\mH^j_i,\ j=1,\ldots,4,$ in order to use the Gaussian off-diagonal decay in $\tilde{q}$ for smoothing the time singularity. Indeed, in this case, for  all $s$ in $[t,T]$ and $y_2$ in $\R^d$, this off-diagonal decay is $\vert y_{2}-m_{t,s}^{2,\xi}(x_1,z_2)\vert$. 

Hence, as we did in Lemma \ref{center2}, we have to re-center the integrands around $(m_{t,s}^{2,\xi}(x_1,\tilde{\xi}_2))_{t<s\leq T}$, for some $\tilde{\xi}_2$ in $\R^{2d}$ with the help of Claim \ref{c1center} and then choose $\tilde{\xi}_2 = z_2$. For notational convenience we write $z_2$ instead of $\tilde{\xi}_2$ as from now. Let 
\begin{eqnarray}
&&\indc \int_{\R^{2d}} \mH_i^j(s,y_1,y_2,\theta_{t,s}^1(\xi),\theta_{t,s}^2(\xi)) D_{z_2}\tilde{q}(t,x_1,z_2;s,y_1,y_2)\d y_1\d y_2\notag\\
&& \quad  := \tilde{\mI}^j_i(s,x_1,z_2) \indc,
\end{eqnarray}
for $j=1,\ldots,4$. Below, we re-center each integrand of $\tilde{\mI}^j_i(s,x_1,z_2)\indc$, $j = 1,\cdots,4,$ and we estimate it.

$\newline$
\emph{Bound of $\tilde{\mI}^1_i(c,x_1,z_2)\indc$.} 
We deduce from assertion \eqref{c1centeresti2} in Claim \ref{c1center} and \eqref{inko} that:
\begin{eqnarray*}
\bigg| \tilde{\mI}^1_i(s,x_1,z_2)\indc \bigg| \leq (s-t)^{-\gamma/2}\int_{\R^{2d}} \Delta^2(\theta_{t,s}(z_2))\phi_i(s,y_1,y_2) D_{z_2}\tilde{q}(t,x_1,z_2;s,y_1,y_2)\d y_1\d y_2  |x_2-z_2|^{\gamma/3} ,
\end{eqnarray*}
Then, by using estimate on $D_{x_2}\tilde{q}$ from Proposition \ref{estfundsol}, regularity of $\phi_i$ under \textbf{(H)} a
\begin{eqnarray}\label{hold_oc_1-2}
\bigg| \tilde{\mI}^1_i(s,x_1,z_2)\indc \bigg| \leq  C\indc (s-t)^{-3(1 - \beta_i^2+\gamma/3)/2}  |x_2-z_2|^{\gamma/3} \label{Final_hold_oc_1-2},
\end{eqnarray}
for all positive $\gamma$ strictly less than $3\beta_i^2-1$.\\

\emph{Bound of $\tilde{\mI}^2_i(s;t,x_1,z_2)\indc$.} We first split this term as:
\begin{eqnarray*}
\tilde{\mI}^2_i(s,x_1,z_2)\indc &=& -\frac{1}{2} \indc\int_{\R^{2d}} {\rm Tr}\left[\left(\Delta^2(m_{t,s}^{t,\xi}(x_1,z_2)) a(s,y)\right)D^2_{y_1}u_i(s,y)\right]D_{z_2}\tilde{q}(t,x_1,z_2;s,y) \d y \\
&&-\frac{1}{2}  \indc\int_{\R^{2d}}  {\rm Tr}\left[ \left(a(s,y_1,m_{t,s}^{2,t,\xi}(x_1,z_2))- a(s,\theta_{t,s}(\xi))\right)D^2_{y_1}u_i(s,y)\right]\\
&& \qquad \times D_{z_2}\tilde{q}(t,x_1,z_2;s,y)\d y .
\end{eqnarray*}
By applying assertion \eqref{c1centeresti2} in Claim \ref{c1center} on the first term in the right hand side, by integrating by parts (see the proof of Lemma \ref{center}) the second term in the right hand side and then by applying assertion \eqref{c1centeresti3} in Claim \ref{c1center} we get:

\begin{eqnarray*}
&&\tilde{\mI}^2_i(s,x_1,z_2)\indc \\
&&= -\frac{1}{2}\int_{t+\epsilon}^T \indc\int_{\R^{2d}} \bigg\{ {\rm Tr}\left[\left(\Delta^2(m_{t,s}^{t,\xi}(x_1,z_2)) a(s,y_1,y_2)\right)D^2_{y_1}u_i(s,y_1,y_2)\right]\\
&&\qquad -\sum_{l=1}^d \bigg(\left[\frac{\p }{\p y_{1l}}a_{l.}(s,y_1,m_{t,s}^{t,\xi}(x_1,z_2))\right]\cdot \left[\Delta^2(m_{t,s}^{t,\xi}(x_1,z_2))D_{y_1}u_i(s,y_1,y_2)\right] \bigg)\bigg\}\\
&& \qquad \hphantom{space}  \times  \left[D_{z_2}\tilde{q}(t,x_1,z_2;s,y_1,y_2)\right]\d y_1 \d y_2 \d s\nonumber\\
&&\quad +\frac{1}{2}\sum_{l=1}^d \int_{t+\epsilon}^T\indc \int_{\R^{2d}} \bigg\{\left[ a_{l.}(s,y_1,m_{t,s}^{2,t,\xi}(x_1,z_2)) - a_{l.}(s,\theta_{t,s}^1(\xi),\theta_{t,s}^2(\xi))\right]\\
&& \qquad \cdot \left[\Delta^2(m_{t,s}^{t,\xi}(x_1,z_2))D_{y_1}u_i(s,y_1,y_2)\right] \times \left[\frac{\p}{\p {y_{1l}}}D_{z_2} \tilde{q}(t,x_1,z_2;s,y_1,y_2)\right] \bigg\}\d y_1 \d y_2 \d s.\nonumber
\end{eqnarray*}
Therefore, from the estimates on the derivatives of $\tilde{q}$ from Proposition \ref{estfundsol}, by using MVT \eqref{labelle1} and the regularity of $a$ from \textbf{(H)} together with estimate \eqref{gaussdecay} and by integrating next in space we obtain that
\begin{eqnarray}\label{hold_oc_2-2}
&&\left|\tilde{\mI}^2_i(s,x_1,z_2)\indc \right|  \notag\\
&&\leq   C' \indc \bigg\{ \bigg( |D^2_{1}u_i||_{\infty}(s-t)^{-\gamma/2}  +||D_{1}D_{2}u_i||_{\infty}  (s-t)^{-\gamma/2} \bigg)|x_2-z_2|^{\gamma/3} \notag\\
&& \quad   +  ||D_{1}D_{2}u_i||_{\infty}   (s-t)^{-1/2} |m_{t,s}^{2,t,\xi}(x_1,z_2)-\theta_{t,s}^{2}(\xi)|   \bigg\},
\end{eqnarray}
by using next the  estimate \eqref{inko} in the first term on the right hand side. This  holds for all $\gamma<2$.\\

\emph{Bound of $\tilde{\mI}^3_i(s,x_1,z_2)\indc$.}  By using the definition of $\Delta^2$ and then assertion \eqref{c1centeresti3} in Claim \ref{c1center}, this term can be centered as follows
\begin{eqnarray*}
\tilde{\mI}^3_i(s,x_1,z_2)\indc\notag &=&- \indc\int_{\R^{2d}} \bigg\{\left[\Delta^2(m_{t,s}^{t,\xi}(x_1,z_2))F_1(s,y_1,y_2)\right] \cdot \left[D_{y_1}u_i(s,y_1,y_2)\right] \\
&& \quad   + \left[F_1(s,y_1,m_{t,s}^{2,t,\xi}(x_1,z_2)) -F_1(s,\theta_{t,s}^1(\xi),\theta_{t,s}^2(\xi))\right]\\
&& \qquad \cdot \left[\Delta^2(m_{t,s}^{2,t,\xi}(x_1,z_2))D_{y_1}u_i(s,y_1,y_2)\right]\bigg\}D_{z_2}\tilde{q}(t,x_1,z_2;s,y_1,y_2)\d y_1\d y_2 .
\end{eqnarray*}
Now, thanks to the regularity of $F_1$ assumed in \textbf{(H)}, if we apply the estimate on the derivative of $\tilde{q}$ from Proposition \ref{estfundsol}, MVT \eqref{labelle1} and estimate \eqref{inko}, we obtain that

\begin{eqnarray}\label{hold_oc_3-2}
&&\left|\tilde{\mI}^3_i(s,x_1,z_2)\indc\right|\\
&& \leq  C'' \indc \bigg\{\bigg(||D_{1}u_i||_{\infty} (s-t)^{-3(1 - \beta_1^2+\gamma/3)/2}   + ||D_{1}D_{2}u_i||_{\infty} (s-t)^{(\beta_1^1- \gamma)/2} \bigg) \notag\\
&& \quad \times |x_2-z_2|^{\gamma/3}+ ||D_{1}D_{2}u_i||_{\infty} |m_{t,s}^{2,t,\xi}(x_1,z_2)-\theta_{t,s}^{2}(\xi)|^{\beta_1^2} \bigg\} \notag
\end{eqnarray}
for all $\gamma < 3\beta_1^2-1$.\\
$\newline$
\emph{Bound of $\tilde{\mI}^4_i(s,x_1,z_2)\indc$.} From representation \eqref{repsol}, by using the definition of $\Delta^2$ and then by centering the term $D_{2}u_i$ thanks to assertion \eqref{c1centeresti3} in Claim \ref{c1center} we can write
\begin{eqnarray*}
&&\tilde{\mI}^4_i(s,x_1,z_2)\indc\notag\\
&&= -\indc\int_{\R^{2d}}\bigg\{ \left[\Delta^2(m_{t,s}^{t,\xi}(x_1,z_2))F_2(s,y_1,y_2)\right] \cdot \left[D_{y_2}u_i(s,y_1,y_2)\right]\\
&& \quad  +  \bigg[F_2(s,y_1,m_{t,s}^{2,t,\xi}(x_1,z_2)) -F_2(s,\theta_{t,s}^1(\xi),\theta_{t,s}^2(\xi)) - D_{1}F_2(s,\theta_{t,s}^1(\xi),\theta_{t,s}^2(\xi))\Delta^1(\theta_{t,s}(\xi))y_1\bigg]\\
&& \qquad  \cdot \left[\Delta^2(m_{t,s}^{2,t,\xi}(x_1,z_2))D_{y_2}u_i(s,y_1,y_2)\right]\bigg\} D_{z_2}\tilde{q}(t,x_1,z_2;s,y_1,y_2)\d y_1\d y_2.
\end{eqnarray*}
By using the regularity of the coefficients from \textbf{(H)}, \eqref{taylorexp} and estimate on $D_{2}\tilde{q}$ from Proposition \ref{estfundsol} and \eqref{inko} we have:

\begin{eqnarray}
&&\Bigg|\tilde{\mI}^4_i(s,x_1,z_2)\indc\Bigg|\label{hold_oc_4-2}\\
&& \leq C'''\indc \Bigg\{  ||D_{2}u_i||_{\infty} (s-t)^{-3(1 - \beta_2^2+\gamma/3)/2} |x_2-z_2|^{\gamma/3} + M(D_{2}u_i,T)\bigg[|m_{t,s}^{2,t,\xi}(x_1,z_2) - \theta_{t,s}^2(\xi)| \notag \\
&& \qquad +  (s-t)^{-1+\alpha^1/2} \bigg(|x_2-z_2|^{\gamma/3}+|x_2-z_2|^{\beta^2_1} + |x_2-z_2|^{\beta^2_2}+ |x_2-z_2|\bigg)\bigg]\Bigg\}  ,\notag
\end{eqnarray}
for all $\gamma < 3\beta_1^2/2-1$.\\

Now, note that from the definition \eqref{meanGauss} of $m$,
\begin{equation*}
\forall s \in [t,T],\quad m_{t,s}^{2,t,x}(x_1,z_2)-\theta^2_{t,s}(x) = z_2-x_2. \label{diffcoef}
\end{equation*}
Hence, by letting $\xi=x$  in \eqref{hold_oc_1-2}, \eqref{hold_oc_2-2}, \eqref{hold_oc_3-2} and \eqref{hold_oc_4-2}  and combining the resulting estimates with \eqref{hudgeinter}, we deduce that there exist a positive constant $C''''$ and a positive number $\delta_{\ref{holdregderivsecdegen}}^2$, depending only on known parameters in \textbf{(H)},  such that:

\begin{eqnarray}
\int_t^T \left| \mP_s(\mathcal{S}^c)\right| \d s
&\leq &  C''''T^{\delta_{\ref{holdregderivsecdegen}}^2}\left( M(D_{2}u_i,T)+ ||D_{2}u_i||_{\infty} + ||D_{1}D_{2}u_i||_{\infty}+ 1\right)\label{la2}\\
&&\quad \times\left(|x_2-z_2|^{\gamma/3}+|x_2-z_2|^{\beta_2^2}+|x_2-z_2|^{\beta_1^2}+|x_2-z_2|\right), \notag
\end{eqnarray}
for all $\gamma < 3\inf (\beta_1^2,\beta_2^2)-1$.\\

\textbf{``H\"{o}lder estimate'' on $D_{2}u_i$.} Finally, by plugging estimates \eqref{la1} and \eqref{la2} in \eqref{sumterms}, we deduce that there exist a positive constant $C$ and a positive number $\delta_{\ref{holdregderivsecdegen}}$, depending only on known parameters in \textbf{(H)}, such that:
\begin{eqnarray*}
\left|D_{2}u_i(t,x_1,x_2)-D_{2}u_i(t,x_1,z_2)\right| &\leq & C T^{\delta_{\ref{holdregderivsecdegen}}} \left(||D_{2}u_i||_{\infty}+||D_{1}D_{2}u_i||_{\infty}+ M(D_{2}u_i,T)+1 \right)\notag\\
&& \times \left( |x_2-z_2|^{\gamma/3}+|x_2-z_2|^{\beta_2^2}+|x_2-z_2|^{\beta_1^2}+|x_2-z_2|\right).
\end{eqnarray*}
A circular argument concludes the proof of Lemma \ref{holdregderivsecdegen}.\qed\\

\emph{Proof of Claim \ref{estifundasolholder}.} Let $(t<s,x,y)$ in $[0,T]^2\times \R^{2d} \times \R^{2d}$, by using MVT and the Gaussian estimate of $D^2_{2}\tilde{q}$ from Proposition \ref{estfundsol} we have:

\begin{eqnarray}
&&\left|(D_{2}\tilde{q})(t,x_1,x_2;s,y_1,y_2)-(D_{2}\tilde{q})(t,x_1,z_2;s,y_1,y_2)\right|\notag\\
&& \qquad \leq \sup_{\rho \in (0,1)}\left|(D^2_{2}\tilde{q})(t,x_1,x_2+\rho(x_2-z_2);s,y_1,y_2)\right| \left|x_2-z_2\right|\notag\\
&& \qquad \leq C' (s-t)^{-3} \sup_{\rho \in (0,1)}\hat{q}_{\bar{c}}(t,x_1,x_2+\rho(x_2-z_2);s,y_1,y_2) \left|x_2-z_2\right| ,\label{rala1}
\end{eqnarray}
where $\bar{c}$ is a positive constant depending only on known parameters in \textbf{(H)}. Note that on $\mathcal{S}$:
\begin{equation}\label{rala2}
\sup_{\rho \in (0,1)}\hat{q}_{\bar{c}}(t,x_1,x_2+\rho(x_2-z_2);s,y_1,y_2) \leq C'' \hat{q}_{c}(t,x_1,x_2;s,y_1,y_2).
\end{equation}
Combining (\ref{rala1}) and (\ref{rala2}), we obtain:
\begin{eqnarray*}
&&\left|(D_{2}\tilde{q})(t,x_1,x_2;s,y_1,y_2)-(D_{2}\tilde{q})(t,x_1,z_2;s,y_1,y_2)\right|\nonumber\\
&& \quad \leq C'''(s-t)^{-3}\hat{q}_{c}(t,x_1,x_2;s,y_1,y_2)\left|x_2-z_2\right|.
\end{eqnarray*}
Rewrite $ |x_2-z_2| = |x_2-z_2|^{1-\gamma/3} |x_2-z_2|^{\gamma/3}$. Since $|x_2-z_2|<(s-t)^{3/2}$ we have $|x_2-z_2| < (s-t)^{3/2-\gamma/2}|x_2-z_2|^{\gamma/3}$ and \eqref{vmtdt} follows.\\

The second assertion follows from the same arguments.\qed

\begin{center}
\textbf{Acknowledgements}
\end{center} 
I thank  Fran\c{c}ois Delarue for his suggestions, large comments and careful reading of the paper as well as the anonymous referees for providing valuable inputs that improved the paper.

\bibliographystyle{plain}
\bibliography{Bib_AIHP716}

\end{document}